\documentclass[12pt]{article}
\usepackage{epsfig}
\usepackage{algorithm,algorithmic}
\usepackage{amsmath}
\usepackage{amssymb}
\usepackage{url}

\def\endproof{\hfill$\square$
\medskip}
\begin{document}
\title{\vspace{-17mm}
Partial Elastic Shape Registration of 3D Surfaces using Dynamic Programming}
\author{\tt\small Javier Bernal$^1$, Jim Lawrence$^{1,2}$\\
$^1${\small \sl National Institute of Standards and Technology,} \\
{\small \sl Gaithersburg, MD 20899, USA} \\
$^2${\small \sl George Mason University,} \\
{\small \sl 4400 University Dr, Fairfax, VA 22030, USA} \\
{\tt\small $\{$javier.bernal,james.lawrence$\}$} \\
{\tt\small @nist.gov \ \ \  lawrence@gmu.edu}}
\date{\ }
\maketitle
\vspace{-13mm}
\begin{abstract}
The computation of the elastic shape registration of two simple surfaces in $3-$dimensional
space and therefore of the elastic shape distance between them has been investigated by
Kurtek, Jermyn, et~al. who have proposed algorithms to carry out this computation.
These algorithms accomplish this by minimizing a distance function between the surfaces in
terms of rotations and reparametrizations of one of the surfaces, the optimization over
reparametrizations using a gradient approach that may produce a local solution. Now minimizing
in terms of rotations and a special subset of the set of reparametrizations, we propose an
algorithm for minimizing the distance function, the optimization over reparametrizations based
on dynamic programming. This approach does not necessarily produce an optimal solution for the
registration and distance problem, but perhaps a solution closer to optimal than the local
solution that an algorithm with a gradient approach for optimizing over the entire set of
reparametrizations may produce. In fact we propose that when computing the elastic shape
registration of two simple surfaces and the elastic shape distance between them with an
algorithm based on a gradient approach for optimizing over the entire set of
reparametrizations, to use as the input initial solution the optimal rotation and
reparametrization computed with our proposed algorithm.
\\[0.2cm]
\end{abstract}
\section{\large Introduction}
\begin{figure}
\centering
\begin{tabular}{cc}
\includegraphics[width=0.4\textwidth]{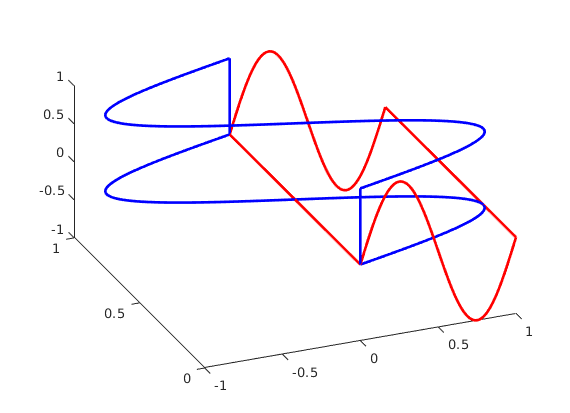}
&
\includegraphics[width=0.4\textwidth]{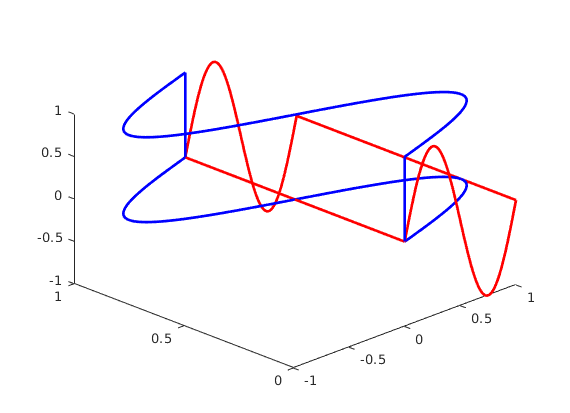}
\end{tabular}
\caption{\label{F:surfaces0}
Two views of the boundaries of the same two surfaces in 3-dimensional space, each of sinusoidal shape.
Their shapes are essentially identical; thus the elastic shape distance between them should be
essentially~zero.
}
\end{figure}

In this paper, we address the problem of computing the elastic shape registration of two simple
surfaces in $3-$dimensional space or equivalently the problem of computing the elastic shape
distance between two such surfaces. Similar work has been carried out by Kurtek, Jermyn et~al.
\cite{kurtek,jermyn}. We do this first through the careful development, independently of
analogous work in \cite{kurtek,jermyn}, of the mathematical framework necessary for the elastic
shape analysis of $3-$dimensional surfaces, which culminates with the definition and justification
of the distance between two such surfaces. This distance, and therefore the registration, is the
result of minimizing a distance function in terms of rotations and reparametrizations of one of
the surfaces. Finally, we propose an algorithm that minimizes the distance function in terms of
rotations and a special subset of the set of reparametrizations, the optimization over
reparametrizations based on Dynamic Programming. Obviously this approach does not necessarily
produce an optimal solution for the registration and distance problem, but perhaps a solution
closer to optimal than the local solution that an algorithm with a gradient approach for
optimizing over the entire set of reparametrizations, such as those proposed in
\cite{kurtek,jermyn}, may produce. In fact we propose that when computing the elastic shape
registration of two simple surfaces and the elastic shape distance between them with an algorithm
based on a gradient approach for optimizing over the entire set of reparametrizations, to use as
the input initial solution the optimal rotation and reparametrization computed with our proposed
algorithm.
\\ \smallskip\\
Given that $S_1$ and $S_2$ are the two surfaces under consideration, we assume they are 
\emph{simple}, that is, we assume elementary regions $D$ and $E$ in the $xy$ plane ($\mathbb{R}^2$)
exist together with one-to-one functions $c_1$ and $c_2$ of class $C^1$,
$c_1:D\rightarrow \mathbb{R}^3$, $c_2:E\rightarrow \mathbb{R}^3$, such that $S_1=c_1(D)$
and $S_2=c_2(E)$. We then say that $c_1$ and $c_2$ \emph{parametrize} or are \emph{parametrizations}
of $S_1$ and $S_2$, respectively,
with domains $D$ and $E$, respectively, and that $S_1$ and $S_2$ are \emph{parametrized surfaces}
relative to $c_1$ and $c_2$, respectively, with domains $D$ and $E$, respectively.
We note that an \emph{elementary region} in the $xy$ plane is one defined by restricting one of
$x$ and $y$ to be between or equal to one of two continuous functions of the remaining variable,
the remaining variable restricted to be in a bounded closed line segment. Actually, for the sake
of simplicity, starting in Section~4 of this paper, we restrict ourselves to exactly one elementary
region, namely $[0,1]\times [0,1]$, the unit square in the $xy$ plane ($\mathbb{R}^2$). Accordingly,
starting in Section~4, we take $D = E = [0,1]\times [0,1]$, and since in practice we can only work
with discretizations of the surfaces $S_1$, $S_2$, given by $c_1$, $c_2$, $D$, $E$ above, we
assume that for positive integers $M$, $N$, not necessarily equal, and
partitions of $[0,1]$, $\{r_i\}_{i=1}^M$, $r_1=0<r_2<\ldots <r_M=1$,
$\{t_j\}_{j=1}^N$, $t_1=0<t_2<\ldots <t_N=1$, not necessarily uniform, $c_1$ and $c_2$ are given
as lists of $M\times N$ points in $S_1$ and $S_2$, respectively, the lists
corresponding to $c_1(r_i,t_j)$ and $c_2(r_i,t_j)$,
$i=1,\ldots,M$, $j=1,\ldots,N$, respectively, and for $k=1,2$, given in the
order $c_k(r_1,t_1)$, $c_k(r_2,t_1)$, $\ldots$, $c_k(r_M,t_1)$, $\ldots$, 
        $c_k(r_1,t_N)$, $c_k(r_2,t_N)$, $\ldots$, $c_k(r_M,t_N)$.
Points $(0,0)$, $(1,0)$, $(1,1)$, $(0,1)$ are the corners of the unit square, and for
$k=1,2$, we can think of $c_k(0,0)$, $c_k(1,0)$, $c_k(1,1)$, $c_k(0,1)$ as the `corners' of the
surface~$S_k$. For the purpose of comparing the shapes of the two surfaces, for each `corner' of $S_2$
we adjust the list of points for~$c_2$ so that the `corner' is the first point in the list, and use
this list together with the list for~$c_1$ to compute a tentative elastic shape distance and
registration between the surfaces. A similar computation is also carried out for the same
`corner' with the adjusted list for $c_2$ in the `reversed' direction (the list for $c_2$ in the
`reversed' direction is given in the order $c_2(r_1,t_1)$, $c_2(r_1,t_2)$, $\ldots$, $c_2(r_1,t_N)$,
$\ldots$, $c_2(r_M,t_1)$, $c_2(r_M,t_2)$, $\ldots$, $c_2(r_M,t_N)$). As $S_2$ has four `corners',
eight tentative elastic shape distances are then obtained and the smallest among them determines the
correct elastic registration of the surfaces. Of course if enough information about the surfaces is
available some of the computations of the tentative elastic shape distances can be avoided and depending
on which take place, $M$ may have to equal $N$, the partitions $\{r_i\}_{i=1}^M$ and $\{t_j\}_{j=1}^N$
may have to be equal, and one or both of them may have to be uniform. For simplicity, in the rest of the
paper, given two simple surfaces $S_1$, $S_2$, as above, we assume the list for $c_2$ suffices as
it is, so that only one tentative elastic shape distance (the correct one) is computed.
\\ \smallskip\\
Being able to compute the elastic shape registration of two surfaces in 3-dimensional space and
the elastic shape distance between them could be useful in studying geological terrains, surfaces of
anatomical objects such as facial surfaces, etc. See Figure~\ref{F:surfaces0} that depicts two such
surfaces (actually their boundaries), each of sinusoidal shape. (Note that in the plots
there, the $x-$, $y-$ and $z-$ axes are not to scale relative to one another).
\\ \smallskip\\
In Section~2 of this paper, we define homeomorphisms and present some known results about them
useful in the context of parametrized simple surfaces in $3-$dimensional space. In particular, we
prove the well-known result that the area of one such surface does not change if its parametrization
is changed. In Section~3, inspired by the definition of the shape function of a parametrized curve
in $d-$dimensional space, $d$~any positive integer, and known results about it, we define the
shape function of a parametrized simple surface in $3-$dimensional space and present some fundamental
results about this function. In Section~4, given two parametrized simple surfaces of unit area in the
form of their shape functions, we associate with them a double integral in terms of rotations of one
of the surfaces, and $C^1$ homeomorphisms with Jacobians of positive determinant, each homeomorphism
corresponding to a reparametrization of the same surface. We then define the elastic shape distance
between the two surfaces as the result of minimizing this double integral with respect to the
aforementioned rotations and homeomorphisms, and justify it accordingly. In Section~5, given two
parametrized simple surfaces of unit area, again in the form of their shape functions, for a fixed
rotation, we describe the computation, based on Dynamic Programming, of a homeomorphism for partially
minimizing the aforementioned double integral, that is, for partially computing the elastic shape
registration of the two surfaces. In Section~6, for a fixed homeomorphism, we describe the computation
of a rotation matrix for approximately minimizing the integral, that is, for approximately computing the
rigid alignment of the two surfaces. In Section~7, we note that the elastic shape distance between
the two surfaces, still in the form of their shape functions, can also be computed in terms of another
double integral that allows for one surface to be reparametrized while the other one is rotated.
We then present a procedure for partially minimizing this other integral, Procedure DP-surface-min,
that alternates computations of optimal homeomorphisms using Dynamic Programming as described in
Section~5, and optimal rotation matrices as described in Section~6. Finally, in Section~8, we present
results obtained with an implementation of our methods.
\section{\large Homeomorphisms and the Area of a Surface}
In this section we present three known results. The first two are about homeomorphisms useful in
the context of parametrized surfaces in $3-$dimensional space, and the third one is about the
computation of the area of one such surface. We note that a \emph{homeomorphism} is a one-to-one
continuous function from a topological space onto another that has a continuous inverse function.
Since simply connected domains are addressed in the first two results that follow, we also note
that a \emph{simply connected domain} is a path-connected domain where one can continuously shrink
any simple closed curve into a point while remaining in the domain. For two-dimensional regions,
a simply connected domain is one without holes in~it. The first result that follows is a standard
result in the field of topology.
\\ \smallskip\\
{\bf Theorem 1:} If $X$ and $Y$ are homeomorphic topological spaces, then $X$ is simply connected
if and only if $Y$ is simply connected.\\ \smallskip\\
{\bf Theorem 2:} Given $D$, a compact simply connected subset of $\mathbb{R}^2$,
and $h:D\rightarrow \mathbb{R}^2$, a homeomorphism, then $h$ maps the boundary of $D$ to exactly
the boundary of~$h(D)$.\\ \smallskip\\
{\bf Proof:} Clearly $h(D)$ is closed as it is a compact subset of $\mathbb{R}^2$, and by Theorem~1
it is simply connected in $\mathbb{R}^2$. Let $p$ be a point in the boundary of~$D$. Then $h$
restricted to $D\setminus \{ p \}$ is a homeomorphism onto~$h(D)\setminus h(p)$. Since
$D\setminus \{ p \}$ is simply connected, it must be that $h(D)\setminus h(p)$ is simply connected
as well so that $h(p)$ cannot be in the interior of~$h(D)$, thus must be in its boundary. On the
other hand, if $q$ is in the boundary of $h(D)$, then through a similar argument since $q$ is in
$h(D)$ ($h(D)$ is closed), it can be shown that $h^{-1}(q)$ is in the boundary of $h^{-1}(h(D)) = D$.
Thus $h$ maps the boundary of $D$ to exactly the boundary of~$h(D)$.
\endproof
\\ \smallskip\\
In what follows, given a surface $S$ in $3-$dimensional space, elementary regions $D$, $E$ in
$\mathbb{R}^2$, and one-to-one functions $c$, $p$ of class $C^1$, $c:D\rightarrow \mathbb{R}^3$,
$p:E\rightarrow \mathbb{R}^3$, $c(D)=S$, $p(E)=S$, so that $c$ and $p$ are parametrizations of $S$
with domain $D$ and $E$, respectively, we say $p$ is a \emph{reparametrization} of $c$ or that $p$
\emph{reparametrizes} $S$ (given as an image of $c$), if $p=c\circ h$ for a $C^1$ homeomorphism $h$
from $E$ onto $D$.  For $(u,v)$ in $D$, writing $c(u,v)=(x(u,v),y(u,v),z(u,v))$, then given a point
$(u_0,v_0)$ in~$D$, the vector tangent to the surface $S$ at $c(u_0,v_0)$ in the $u$ direction is
given by
\[ \frac{\partial c}{\partial u} (u_0,v_0) =
(\frac{\partial x}{\partial u} (u_0,v_0),
\frac{\partial y}{\partial u} (u_0,v_0),
\frac{\partial z}{\partial u} (u_0,v_0)),
\]
and in the $v$ direction by
\[ \frac{\partial c}{\partial v} (u_0,v_0) =
(\frac{\partial x}{\partial v} (u_0,v_0),
\frac{\partial y}{\partial v} (u_0,v_0),
\frac{\partial z}{\partial v} (u_0,v_0)).
\]
We say the surface $S$ is \emph{regular} (relative to the parametrization $c$) if at every point
$c(u_0,v_0)$ in $S$ the cross product
$\frac{\partial c}{\partial u} (u_0,v_0) \times \frac{\partial c}{\partial v} (u_0,v_0)$
is nonzero. We note that if $S$ is regular, then at every point $c(u_0,v_0)$ in $S$, 
$\frac{\partial c}{\partial u} (u_0,v_0) \times \frac{\partial c}{\partial v} (u_0,v_0)$
is a nonzero vector normal to $S$ at $c(u_0,v_0)$.\\ \smallskip\\
With $c$, $D$, $S$, $\frac{\partial c}{\partial u} (u_0,v_0)$, $\frac{\partial c}{\partial v} (u_0,v_0)$
as above, $S$ regular (relative to $c$), the surface area $A(S)$ of the parametrized surface $S$ is
given by \[ A(S) = \int\int_D ||
\frac{\partial c}{\partial u} (u,v) \times \frac{\partial c}{\partial v} (u,v)||\,du\,dv \]
where $||\cdot||$ is the 3-dimensional Euclidean norm.
\\ \smallskip\\
With $c$, $D$, $p$, $E$, $S$, $\frac{\partial c}{\partial u} (u_0,v_0)$,
$\frac{\partial c}{\partial v} (u_0,v_0)$, $h$ as above so that $p$ is also a parametrization of $S$
with domain $E$, and $p$ is a reparametrization of $c$, $p=c\circ h$, the result that follows shows
the surface area $A(S)$ of $S$ does not change if it is computed with the parametrization $p$ of $S$
with domain $E$ instead of the parametrization $c$ of $S$ with domain~$D$. For $(r,t)$ in $E$,
writing $h(r,t)=(u(r,t),v(r,t))$, and letting $\frac{\partial (u,v)}{\partial (r,t)}$
be the determinant of the Jacobian of~$h$, $\frac{\partial (u,v)}{\partial (r,t)}$
is assumed to be nonzero on~$E$.
Finally, for $(r,t)$ in $E$, writing $p(r,t)=(\hat{x}(r,t),\hat{y}(r,t),\hat{z}(r,t))$, then given a
point $(r_0,t_0)$ in~$E$, the
vector tangent to the surface $S$ at $p(r_0,t_0)$ in the $r$ direction is given by
\[ \frac{\partial p}{\partial r} (r_0,t_0) =
(\frac{\partial \hat{x}}{\partial r} (r_0,t_0),
\frac{\partial \hat{y}}{\partial r} (r_0,t_0),
\frac{\partial \hat{z}}{\partial r} (r_0,t_0)),
\]
and in the $t$ direction by
\[ \frac{\partial p}{\partial t} (r_0,t_0) =
(\frac{\partial \hat{x}}{\partial t} (r_0,t_0),
\frac{\partial \hat{y}}{\partial t} (r_0,t_0),
\frac{\partial \hat{z}}{\partial t} (r_0,t_0)).
\]
\\ \smallskip\\
{\bf Theorem 3:} Given $c$, $D$, $p$, $E$, $S$,
$\frac{\partial c}{\partial u} (u_0,v_0)$, $\frac{\partial c}{\partial v} (u_0,v_0)$,
$\frac{\partial p}{\partial r} (r_0,t_0)$, $\frac{\partial p}{\partial t} (r_0,t_0)$,
$h$, $\frac{\partial (u,v)}{\partial (r,t)}$ as above,
then $S$ is regular relative to $p$ and
\[
\int\int_E || \frac{\partial p}{\partial r} (r,t) \times \frac{\partial p}{\partial t} (r,t)||\,dr\,dt =
\int\int_D || \frac{\partial c}{\partial u} (u,v) \times \frac{\partial c}{\partial v} (u,v)||\,du\,dv.
\]
\smallskip \\
{\bf Proof:} With $(u(r,t),v(r,t))=h(r,t)$, then $$p(r,t) = c(h(r,t)) = c(u(r,t),v(r,t)),$$ so that 
\[
\frac{\partial p}{\partial r}=
\frac{\partial c}{\partial u} \frac{\partial u}{\partial r} +
\frac{\partial c}{\partial v} \frac{\partial v}{\partial r},
\ \ \ \ \ \frac{\partial p}{\partial t}=
\frac{\partial c}{\partial u} \frac{\partial u}{\partial t} +
\frac{\partial c}{\partial v} \frac{\partial v}{\partial t}.
\] 
Thus
\begin{eqnarray*}
\frac{\partial p}{\partial r} \times \frac{\partial p}{\partial t} &=&
(\frac{\partial c}{\partial u} \frac{\partial u}{\partial r} +
 \frac{\partial c}{\partial v} \frac{\partial v}{\partial r}) \times
(\frac{\partial c}{\partial u} \frac{\partial u}{\partial t} +
 \frac{\partial c}{\partial v} \frac{\partial v}{\partial t})\\
&=&(\frac{\partial c}{\partial u} \times \frac{\partial c}{\partial v})
(\frac{\partial u}{\partial r} \frac{\partial v}{\partial t}) +
(\frac{\partial c}{\partial v} \times \frac{\partial c}{\partial u})
(\frac{\partial v}{\partial r} \frac{\partial u}{\partial t})\\
&=&(\frac{\partial c}{\partial u} \times \frac{\partial c}{\partial v})
(\frac{\partial u}{\partial r} \frac{\partial v}{\partial t} -
 \frac{\partial v}{\partial r} \frac{\partial u}{\partial t})\\
&=&(\frac{\partial c}{\partial u} \times \frac{\partial c}{\partial v})
\;\frac{\partial (u,v)}{\partial (r,t)}
\end{eqnarray*}
so that $S$ is regular relative to $p$ since
both $\frac{\partial c}{\partial u} \times \frac{\partial c}{\partial v}$ and
$\frac{\partial (u,v)}{\partial (r,t)}$ are nonzero on~$E$, and
\begin{eqnarray*}
\int\int_E || \frac{\partial p}{\partial r} \times \frac{\partial p}{\partial t} ||\,dr\,dt
&=& \int\int_E || (\frac{\partial c}{\partial u} \times \frac{\partial c}{\partial v})
\;\frac{\partial (u,v)}{\partial (r,t)} ||\,dr\,dt\\
&=& \int\int_E || \frac{\partial c}{\partial u} \times \frac{\partial c}{\partial v} ||
\;|\frac{\partial (u,v)}{\partial (r,t)}|\,dr\,dt\\
&=& \int\int_D || \frac{\partial c}{\partial u} \times \frac{\partial c}{\partial v} ||\,du\,dv
\end{eqnarray*}
by the change of variables formula. \endproof
\section{\large The Shape Function of a Parametrized Surface}
In this section we define the shape function of a parametrized surface in 3-dimensional space
and present some fundamental results about this function. A similar definition and similar
results have been presented in \cite{bernal2,bernal4,joshi,srivastava,srivastava2} in the
context of the shape function of a parametrized curve in $d-$dimensional space, $d$ any
positive integer. Accordingly, in \cite{bernal2,bernal4,joshi,srivastava,srivastava2}, given
$\beta:[0,1]\rightarrow \mathbb{R}^d$ of class $C^1$, a parametrization of a curve
in~$\mathbb{R}^d$, the shape function $q$ of $\beta$, i.e., the shape function $q$ of the
curve that $\beta$ parametrizes relative to $\beta$, $q:[0,1]\rightarrow \mathbb{R}^d$, is
defined by $q(t)=\dot{\beta}(t)/\sqrt{||\dot{\beta}(t)||}$, $t\in [0,1]$ ($d-$dimensional~0 if
$\dot{\beta}(t)$ equals $d-$dimensional~0).
It follows then that $q$ is square integrable as
\[ \int_0^1 ||q(t)||^2 dt=\int_0^1 ||\dot{\beta}(t)/\sqrt{||\dot{\beta}(t)||}\,||^2 dt=
\int_0^1 ||\dot{\beta}(t)|| dt \]
which is the length of the curve that $\beta$ parametrizes, where $||\cdot||$ is the
$d-$dimensional Euclidean norm.
Note that in what follows we ignore the usual definition of a diffeomorphism and refer to
homeomorphisms on~$[0,1]$ as diffeomorphisms in order
to distinguish them from homeomorphisms on elementary regions in the plane.
Again with $q$ the shape function of $\beta$ and $\Gamma$ the set
of $C^1$ orientation-preserving diffeomorphisms of $[0,1]$ so that for $\gamma\in\Gamma$ then
$\dot{\gamma}\geq 0$ on $[0,1]$, it then follows that for $\gamma\in\Gamma$ the shape function of
the reparametrization $\beta\circ\gamma$ of $\beta$ is
$(q,\gamma)=(q\circ\gamma)\sqrt{\dot{\gamma}}$.  With $||q||_2 = (\int_0^1 ||q(t)||^2dt)^{1/2}$,
we also note that given $\beta_1$, $\beta_2:[0,1]\rightarrow \mathbb{R}^d$ of class $C^1$,
parametrizations of curves in $\mathbb{R}^d$ with shape functions $q_1$, $q_2$, respectively,
then $||(q_1,\gamma) - (q_2,\gamma)||_2=||q_1-q_2||_2$ for any $\gamma\in\Gamma$, and from this,
with $\Gamma_0=\{\gamma\in\Gamma,\ \dot{\gamma}>0$ on $[0,1]\}$, it has been demonstrated
\cite{bernal2,srivastava} that ignoring rotations, the number
inf$_{\gamma\in\Gamma_0}||q_1-(q_2,\gamma)||_2$ can then be used as a well-defined distance
between the two curves that $\beta_1$,~$\beta_2$ parametrize, $\beta_1$ and $\beta_2$ both
normalized to parametrize curves of length~1.
\\ \smallskip \\
With $c$, $D$, $S$, $\frac{\partial c}{\partial u}$, $\frac{\partial c}{\partial v}$ as in the
previous section, $S$ regular (relative to $c$), following the idea of the definition of the shape
function of a parametrized curve in $d-$dimensional space as described above,
we define the shape function $q$ of the parametrization $c$ of $S$ with domain $D$, i.e., the
shape function $q$ of the surface $S$ relative to its parametrization $c$ with domain~$D$,
$q:D\rightarrow \mathbb{R}^3$, so that \[ \int\int_D ||q(u,v)||^2 du\,dv =
\int\int_D || \frac{\partial c}{\partial u} \times \frac{\partial c}{\partial v} ||\,du\,dv
\]
which is the surface area of~$S$. This is easily seen to be indeed the case if we define the
shape function $q$ of $c$ on $D$ by
\[ q=(\frac{\partial c}{\partial u} \times \frac{\partial c}{\partial v})/
\sqrt{|| \frac{\partial c}{\partial u} \times \frac{\partial c}{\partial v}||}. \]
\smallskip\\
We do define $q$ this way and note that this definition of the shape function of a surface relative
to a parametrization of the surface, is slightly different from the one in~\cite{kurtek} but
similar to the one in \cite{jermyn}. We also note that if we allow
$\frac{\partial c}{\partial u} \times \frac{\partial c}{\partial v}$ to be $3-$dimensional zero 
at certain points, then $q$ is defined to be $3-$dimensional zero at those points.
\\ \smallskip\\
With $c$, $q$, $D$, $S$ as above, the following result, similar to the one mentioned above in the
context of the shape function of the parametrization of a curve in $d-$dimensional space, shows
how to compute the shape function of a reparametrization of~$c$ from the shape function $q$ of~$c$.
Here $p$ is the reparametrization of $c$, i.e., for an elementary region $E$ in $\mathbb{R}^2$,
$p$ is a parametrization of $S$ with domain $E$, and $p=c\circ h$ for a $C^1$ homeomorphism $h$
from $E$ onto $D$. Assuming $\frac{\partial (u,v)}{\partial (r,t)}\geq 0$ on~$E$,
$\frac{\partial (u,v)}{\partial (r,t)}$ the determinant of the Jacobian of~$h$, we define a
function on $E$ into $\mathbb{R}^3$, which we denote by $(q,h)$, as follows:
\[ (q,h)\equiv (q\circ h)\sqrt{\frac{\partial (u,v)}{\partial (r,t)}}. \]
\\ \smallskip\\
{\bf Theorem 4:} Given $c$, $q$, $D$, $p$, $h$, $E$, $S$, $\frac{\partial (u,v)}{\partial (r,t)}$
as above, the shape function on $E$ of the reparametrization $p=c\circ h$ of $c$ is then~$(q,h)$.
\smallskip\\
{\bf Proof:} With $(u(r,t),v(r,t))=h(r,t)$, so that $$p(r,t) = c(h(r,t)) = c(u(r,t),v(r,t)),$$
then on $E$, as established in the proof of Theorem~3, we have
\[ \frac{\partial p}{\partial r} \times \frac{\partial p}{\partial t} =
(\frac{\partial c}{\partial u} \times \frac{\partial c}{\partial v})
\;\frac{\partial (u,v)}{\partial (r,t)}. \]
Thus, if $\hat{q}$ is the shape function of $p$ on $E$, from the definition of a shape function
it must then be that
\begin{eqnarray*}
\hat{q} &=&(\frac{\partial p}{\partial r} \times \frac{\partial p}{\partial t})/
\sqrt{|| \frac{\partial p}{\partial r} \times \frac{\partial p}{\partial t}||}\\
&=& (\frac{\partial c}{\partial u} \times \frac{\partial c}{\partial v})
\;\frac{\partial (u,v)}{\partial (r,t)}/
\sqrt{||(\frac{\partial c}{\partial u} \times \frac{\partial c}{\partial v})
\;\frac{\partial (u,v)}{\partial (r,t)}||}\\
&=& \big(\;(\frac{\partial c}{\partial u} \times \frac{\partial c}{\partial v})/
\sqrt{||\frac{\partial c}{\partial u} \times \frac{\partial c}{\partial v}||}\;\big)
\;\sqrt{\frac{\partial (u,v)}{\partial (r,t)}}\\
&=&(q\circ h)\sqrt{\frac{\partial (u,v)}{\partial (r,t)}} = (q,h).
\end{eqnarray*}
\endproof
\smallskip\\
Other results about shape functions of parametrized surfaces, similar to results about
shape functions of parametrized curves in $d-$dimensional space \cite{bernal2,srivastava},
can be developed.  Given $c$, $q$, $p$, $h$, $D$, $E$, $S$, $\frac{\partial (u,v)}{\partial (r,t)}$
as above, $q$ the shape function of $c$,
$p=c\circ h$, $h$ a $C^1$ homeomorphism from $E$ onto $D$,
$\frac{\partial (u,v)}{\partial (r,t)}$ the determinant of the Jacobian of~$h$, assuming now
$\frac{\partial (u,v)}{\partial (r,t)}>0$ on~$E$, so that
$\frac{\partial (r,t)}{\partial (u,v)}(u,v)=(\frac{\partial (u,v)}{\partial (r,t)}(r,t))^{-1}$,
$\frac{\partial (r,t)}{\partial (u,v)}$ the determinant of the Jacobian of~$h^{-1}$ on~$D$,
then with $(q,h)$ as defined above,
one such result is that $((q,h),h^{-1}) = q$ on~$D$. This result together with the theorem that
follows are of importance in the next section for justifying the definition of the distance between
surfaces in a manner similar to the way the definition of the distance between curves in
$d-$dimensional space is justified~\cite{bernal2,srivastava}. The theorem shows homeomorphisms act
by isometries on shape functions of parametrized surfaces.
\\ \smallskip\\
{\bf Theorem 5:} Given $D$, $E$, $h$, $\frac{\partial (u,v)}{\partial (r,t)}$ as above,
$\frac{\partial (u,v)}{\partial (r,t)}\geq 0$ on~$E$;
$S_1$, $S_2$ surfaces, $c_1$, $c_2$ parametrizations of $S_1$, $S_2$, respectively, both with domain~$D$;
$p_1$, $p_2$ parametrizations of $S_1$, $S_2$, respectively, both with domain~$E$; $p_1$, $p_2$
reparametrizations of $c_1$, $c_2$, respectively, $p_1=c_1\circ h$, $p_2=c_2\circ h$; $q_1$, $q_2$,
$\hat{q}_1$,~$\hat{q}_2$ the shape functions of $c_1$, $c_2$, $p_1$, $p_2$, respectively, then
\begin{eqnarray*}
||\hat{q}_1-\hat{q}_2||_{2,E} &\equiv& \big(\int\int_E ||\hat{q}_1-\hat{q}_2)||^2dr\,dt\big)^{1/2}\\
&=& \big(\int\int_D ||q_1-q_2||^2du\,dv\big)^{1/2}\\
&\equiv& ||q_1-q_2||_{2,D}.
\end{eqnarray*}
\smallskip\\
{\bf Proof:} From Theorem 4, $\hat{q}_1 = (q_1,h)$, $\hat{q}_2 = (q_2,h)$, thus
\begin{eqnarray*}
||\hat{q}_1-\hat{q}_2||_{2,E}^2 &=& \int\int_E ||\hat{q}_1-\hat{q}_2||^2dr\,dt\\
&=& \int\int_E ||(q_1,h)-(q_2,h)||^2dr\,dt\\
&=& \int\int_E ||(q_1\circ h)\sqrt{\frac{\partial (u,v)}{\partial (r,t)}}
-(q_2\circ h)\sqrt{\frac{\partial (u,v)}{\partial (r,t)}}||^2dr\,dt\\
&=& \int\int_E ||(q_1\circ h) -(q_2\circ h)||^2\;\frac{\partial (u,v)}{\partial (r,t)}dr\,dt\\
&=& \int\int_D ||q_1-q_2||^2du\,dv\\
&=& ||q_1-q_2||_{2,D}^2
\end{eqnarray*}
by the change of variables formula.
\endproof
\section{\large The Elastic Shape Distance between Surfaces}
In this section we define and justify the elastic shape distance between two surfaces of unit area.
This is done at first in terms of $C^1$ homeomorphisms with Jacobians of
positive determinant (each homeomorphism defines a reparametrization of one of the surfaces),
and later in terms of rotations as well. Given that $S_1$ and $S_2$ are the
two surfaces, we assume they are simple and are parametrized by functions with the same domain,
i.e., an elementary region $D$ in the $xy$ plane exists together with
parametrizations $c_1$ and $c_2$ with domain $D$ of $S_1$ and $S_2$, respectively,
$c_1:D\rightarrow \mathbb{R}^3$, $c_2:D\rightarrow \mathbb{R}^3$, $S_1=c_1(D)$, $S_2=c_2(D)$.
Letting $\Sigma_0$ be the set of all $C^1$ homeomorphisms $h$, from $D$ onto $D$,
with $\frac{\partial (u,v)}{\partial (r,t)}>0$ on~$D$,
$\frac{\partial (u,v)}{\partial (r,t)}$ the determinant of the Jacobian of~$h$,
given that $q_1$ and $q_2$ are, respectively, the shape functions of $c_1$ and $c_2$, then
using arguments similar to arguments for justifying the definition of the distance between
curves in $d-$dimensional space found in \cite{bernal2,srivastava}, ignoring rotations, it can be
demonstrated that the number
\[ \mathrm{inf}_{h\in\Sigma_0}||q_1-(q_2,h)||_{2,D} =
\mathrm{inf}_{h\in\Sigma_0}\big(\int\int_D ||q_1-(q_2,h)||^2dr\,dt\big)^{1/2} = \]
\[ \mathrm{inf}_{h\in\Sigma_0}\big(\int\int_D ||q_1
-(q_2\circ h)\sqrt{\frac{\partial (u,v)}{\partial (r,t)}}||^2dr\,dt\big)^{1/2} \]
can be used as a well-defined distance between the surfaces $S_1$ and~$S_2$,
$c_1$ and $c_2$ both normalized to parametrize surfaces of area equal to~1.
Note that the arguments for justifying this definition of the distance between the two
surfaces are in part based on Theorem~$5$ in the previous section and the result described
in the paragraph preceding Theorem~5.
\\ \smallskip\\
That rotations as well act by isometries on shape functions of parametrized surfaces is
justified as follows. With $q_1$, $q_2$ as above, assuming $R$ is a $3-$dimensional rotation
matrix, i.e., $R\in SO(3)$, $SO(3)$ the group of $3\times 3$ orthogonal matrices of determinant
equal to~1, then because $R$ is orthogonal, it follows easily that
\begin{eqnarray*}
||Rq_1-Rq_2||_{2,D} &=& \big(\int\int_D ||Rq_1-Rq_2)||^2du\,dv\big)^{1/2}\\
&=& \big(\int\int_D ||q_1-q_2||^2du\,dv\big)^{1/2}\\
&=& ||q_1-q_2||_{2,D}.
\end{eqnarray*}
Also as established in \cite{srivastava} for shape functions of parametrized curves in
$d-$dimensional space, it follows by similar arguments that given $h\in\Sigma_0$, $R\in SO(3)$,
$D$ an elementary region in the $xy$ plane, $q$ a shape function of a surface parametrized by
a function $c$ from $D$ into $\mathbb{R}^3$, then $(Rq,h) = R(q,h)$. That is, the actions
on shape functions of homeomorphisms in $\Sigma_0$ and matrices in $SO(3)$ commute.
For the sake of completeness we actually present the details of the justification of this fact in
what follows. However, for this purpose, we first present a well-known formula about rotations and
cross products of vectors in $\mathbb{R}^3$ together with its justification, again for the
sake of completeness.\\ \smallskip\\
{\bf Lemma:} Given vectors $x$, $y$ in $\mathbb{R}^3$, $R$ in $SO(3)$,
then $R(x\times y) = Rx\times Ry$.\\ \smallskip\\
{\bf Proof:} Here given $a$, $b$, $c$ in $\mathbb{R}^3$, we use the identity
$a\cdot (b\times c) = \det [a\;b\;c]$, where $\cdot$ and $\det$ denote the inner product and
determinant operations, respectively. In addition, for $j=1,2,3$, we let $e_j$ be the $j^{th}$ unit
vector in $\mathbb{R}^3$, and given $w$ in $\mathbb{R}^3$, we let $(w)_j$ denote the $j^{th}$
coordinate of~$w$. Since $\det R=1$ and $RR^T$ equals the identity matrix, we then have
for $j=1,2,3$,
\begin{eqnarray*}
(R(x\times y))_j &=& e_j\cdot R(x\times y) = e_j^T R(x\times y) = (R^Te_j)^T (x\times y)\\
&=& R^Te_j\cdot (x\times y) = \det [R^Te_j\ \;x\ \;y]\\
&=& \det R\;\det [R^Te_j\ \;x\ \;y] = \det R\;[R^Te_j\ \;x\ \;y]\\
&=& \det [RR^Te_j\ \;Rx\ \;Ry] = \det [e_j\ \;Rx\ \;Ry]\\
&=& e_j\cdot (Rx\times Ry) = (Rx\times Ry)_j.
\end{eqnarray*}
Thus $R(x\times y) = Rx\times Ry$. \endproof
\\ \smallskip\\
With $h$, $R$, $q$, $c$ as above, in order to show $(Rq,h)=R(q,h)$, we first show that the shape
function of $Rc$ on $D$, say $\hat{q}$, is $Rq$. From the definition of a shape function and the
lemma then
\begin{eqnarray*}
\hat{q}&=&(\frac{\partial Rc}{\partial u} \times \frac{\partial Rc}{\partial v})/
\sqrt{|| \frac{\partial Rc}{\partial u} \times \frac{\partial Rc}{\partial v}||}
=(R\frac{\partial c}{\partial u} \times R\frac{\partial c}{\partial v})/
\sqrt{|| R\frac{\partial c}{\partial u} \times R\frac{\partial c}{\partial v}||}\\
&=&R(\frac{\partial c}{\partial u} \times \frac{\partial c}{\partial v})/
\sqrt{|| R(\frac{\partial c}{\partial u} \times \frac{\partial c}{\partial v})||}
=R(\frac{\partial c}{\partial u} \times \frac{\partial c}{\partial v})/
\sqrt{|| \frac{\partial c}{\partial u} \times \frac{\partial c}{\partial v}||}\\
&=& Rq.
\end{eqnarray*}
From Theorem 4 and what we just proved, it follows that the shape function of $Rc(h)$ is then
$(Rq,h)$. On the other hand, again by Theorem~4, the shape function of $c(h)$ is $(q,h)$ so that
again by what we just proved the shape function of $R(c(h))$ must be $R(q,h)$. Since $Rc(h)$ and
$R(c(h))$ are the same function, then it must be that their shape functions are the same, i.e.,
$(Rq,h)=R(q,h)$. \\ \smallskip\\
Based in part on the observations above about rotation matrices and homeomorphisms, in a manner
similar to what is done in \cite{jermyn,kurtek}, given $S_1$, $S_2$, $c_1$, $c_2$, $q_1$, $q_2$
as above, with inf short for infimum, it can be demonstrated that the number \smallskip
\[ \mathrm{inf}_{R\in SO(3),h\in\Sigma_0}||q_1-R(q_2,h)||_{2,D} = \]
\[ \mathrm{inf}_{R\in SO(3),h\in\Sigma_0}\big(\int\int_D ||q_1-R(q_2,h)||^2dr\,dt\big)^{1/2} = \]
\[ \mathrm{inf}_{R\in SO(3),h\in\Sigma_0}\big(\int\int_D ||q_1
-R(q_2\circ h)\sqrt{\frac{\partial (u,v)}{\partial (r,t)}}||^2dr\,dt\big)^{1/2} \]
\smallskip\\
can be used as a well-defined distance between the surfaces $S_1$ and~$S_2$,
where again $\frac{\partial (u,v)}{\partial (r,t)}$ is the determinant of the Jacobian of~$h$,
and $c_1$ and $c_2$ are both normalized to parametrize surfaces of area equal to~1.
Thus, denoting $\mathrm{inf}_{R\in SO(3),h\in\Sigma_0}||q_1-R(q_2,h)||_{2,D}$ by dist$(S_1,S_2)$,
and restricting ourselves to the simpler region $D=[0,1]\times [0,1]$, then by Fubini's theorem,
we note,
\begin{eqnarray*}
\mathrm{dist}(S_1,S_2)&=&
\mathrm{inf}_{R\in SO(3),h\in\Sigma_0}\big(\int\int_D ||q_1-R(q_2,h)||^2dr\,dt\big)^{1/2}\\
&=& \mathrm{inf}_{R\in SO(3),h\in\Sigma_0}\big(\int_0^1\int_0^1 ||q_1-R(q_2,h)||^2dr\,dt\big)^{1/2}
\end{eqnarray*}
which we use in the next section.
\section{\large Computation of Homeomorphism for Partial Registration of Surfaces
using Dynamic Programming}
In this section, ignoring rotations, we describe the computation, based on Dynamic Programming,
of a homeomorphism for the partial elastic shape registration of two simple surfaces of
unit area in $3-$dimensional space. Given that $S_1$ and $S_2$ are the two surfaces, with
$D=[0,1]\times [0,1]$, we assume accordingly that one-to-one functions $c_1$ and $c_2$ exist of
class $C^1$, $c_1:D\rightarrow \mathbb{R}^3$, $c_2:D\rightarrow \mathbb{R}^3$, such that
$S_1=c_1(D)$ and $S_2=c_2(D)$. That is, $c_1$ and $c_2$ parametrize or are parametrizations
of $S_1$ and $S_2$, respectively. Given that $q_1$ and $q_2$ are, respectively, the shape
functions of $c_1$ and $c_2$, then we hope to minimize
\begin{eqnarray*}
 \int\int_D ||q_1-(q_2,h)||^2dr\,dt &=&
 \int\int_D ||q_1 -(q_2\circ h)\sqrt{\frac{\partial (u,v)}{\partial (r,t)}}||^2dr\,dt\\
 &=& \int_0^1\int_0^1 ||q_1 -(q_2\circ h)\sqrt{\frac{\partial (u,v)}{\partial (r,t)}}||^2dr\,dt
\end{eqnarray*}
with respect to $h$ in $\Sigma_0$, where $\Sigma_0$ is the set of all $C^1$ homeomorphisms $h$
from $D$ onto itself, with $\frac{\partial (u,v)}{\partial (r,t)}>0$ on~$D$,
$\frac{\partial (u,v)}{\partial (r,t)}$ the determinant of the Jacobian of~$h$.
We note, this minimization is usually carried out with an algorithm that uses a gradient approach
for the optimization over reparametrizations, i.e., over homeomorphisms $h$ in $\Sigma_0$, that may
produce a local solution \cite{kurtek,jermyn}. In this paper we have opted to carry out the minimization
with respect to $h$ in a special subset of $\Sigma_0$ that allows for the use of Dynamic Programming.
We denote this subset of $\Sigma_0$ by $\Sigma_1$, $h$ in $\Sigma_1$ satisfying that $h\in \Sigma_0$
and for $(r,t)$ in $D$, if $h(r,t)=(\hat{r},\hat{t})$
then it must be that $\hat{t}=t$. In addition, if $h \in \Sigma_1$, we assume for any $t$ in $[0,1]$
that $h(0,t)=(0,t)$ and $h(1,t)=(1,t)$. That a minimization over $\Sigma_1$ allows for the use of
Dynamic Programming will become evident below.
Note, from Theorem~2, with $\partial D$ the boundary of~$D$, for any homeomorphism $h$ from $D$
onto itself, not necessarily in $\Sigma_0$ or $\Sigma_1$, it must be that~$h(\partial D) = \partial D$.
\\ \smallskip\\
Since in practice we can only work with a discretized version of the problem, for our purposes
we assume the situation is as follows: for positive integers $M$, $N$, not necessarily equal, and
partitions of $[0,1]$, $\{r_i\}_{i=1}^M$, $r_1=0<r_2<\ldots <r_M=1$,
$\{t_j\}_{j=1}^N$, $t_1=0<t_2<\ldots <t_N=1$, not necessarily uniform, $c_1$ and $c_2$ are given
as lists of $M\times N$ points in the surfaces $S_1$ and $S_2$, respectively, the lists for
$c_1$ and $c_2$ corresponding to $c_1(r_i,t_j)$ and $c_2(r_i,t_j)$,
$i=1,\ldots,M$, $j=1,\ldots,N$, respectively; for $k=1,2$, the list for $c_k$ given in the following
order: $c_k(r_1,t_1)$, $c_k(r_2,t_1)$, $\ldots$, $c_k(r_M,t_1)$, $\ldots$,
        $c_k(r_1,t_N)$, $c_k(r_2,t_N)$, $\ldots$, $c_k(r_M,t_N)$.
\\ \smallskip\\
Computing $\frac{\partial c_1}{\partial r}(r_i,t_j)$, $\frac{\partial c_1}{\partial t}(r_i,t_j)$,
$\frac{\partial c_2}{\partial r}(r_i,t_j)$, $\frac{\partial c_2}{\partial t}(r_i,t_j)$
with centered finite differences from $c_1(r_i,t_j)$ and $c_2(r_i,t_j)$, for $i=1,\ldots,M$,
$j=1,\ldots,N$, we can then approximately compute for $i=1,\ldots,M$, $j=1,\ldots,N$, $k=1,2$,
\[ q_k(r_i,t_j)=\big( \;(\frac{\partial c_k}{\partial r} \times \frac{\partial c_k}{\partial t})/
\sqrt{|| \frac{\partial c_k}{\partial r} \times \frac{\partial c_k}{\partial t}||}\;\big)
(r_i,t_j), \] ($3-$dimensional zero if
$(\frac{\partial c_k}{\partial r} \times \frac{\partial c_k}{\partial t})(r_i,t_j)$ equals
$3-$dimensional zero).
\\ \smallskip\\
Finally we note that if $h\in\Sigma_1$, then $h(r,t)=(u(r,t),v(r,t))=(u(r,t),t)$ for $(r,t)$ in
$D$, so that $\frac{\partial v}{\partial r}(r,t) = 0$ and $\frac{\partial v}{\partial t}(r,t) = 1$,
and therefore $\frac{\partial (u,v)}{\partial (r,t)}(r,t)=\frac{\partial u}{\partial r}(r,t)$ for
$(r,t)$ in $D$. 
\\ \smallskip\\
Given $h$ in $\Sigma_1$ and an integer $j$, $1\leq j \leq N$, next we discretize the integral
 \[ \int_0^1 ||q_1(r,t_j) -\big(\; (q_2\circ h)\sqrt{\frac{\partial (u,v)}{\partial (r,t)}}\;\big)
(r,t_j)||^2dr. \]
For this purpose, we define $q_{1j}(r_i)$, $q_{2j}(r_i)$ in $\mathbb{R}^3$ for $i=1,\ldots,M$, by
\[ q_{1j}(r_i)=q_1(r_i,t_j),\  q_{2j}(r_i)=q_2(r_i,t_j), \]
and define as well a diffeomorphism $h_j$ from $[0,1]$ onto $[0,1]$ by
\[ h_j(r)=u(r,t_j),\ r\in [0,1]. \]
Note, $h_j$ is indeed a diffeomorphism as clearly $h_j(0)=0$, $h_j(1)=1$, and for $r\in [0,1]$,
$h_j'(r) = \frac{dh_j}{dr}(r) = \frac{\partial u}{\partial r}(r,t_j)=
\frac{\partial (u,v)}{\partial (r,t)}(r,t_j)>0$.\\
For $i=1,\ldots,M$, we can then compute $h_j(r_i) = u(r_i,t_j)$ so that $h_j(r_1)=h_j(0)=0$ and
$h_j(r_M)=h_j(1)=1$. In addition, for $i=1,\ldots,M-1$, we compute $\Delta r_i = r_{i+1}-r_i$,
approximately compute $h_j'(r_i) = (h_j(r_{i+1})-h_j(r_i))/\Delta r_i$, set $h_j'(r_M)=h_j'(r_1)$,
and by interpolating $q_{2j}(r_i)$, $i=1,\ldots,M$, by a cubic spline, for $i=1,\ldots,M$, we can
approximately compute $q_{2j}(h_j(r_i))$, which in turn is an approximation of
$(q_2\circ h)(r_i,t_j)$ as $(q_2\circ h)(r_i,t_j) = q_2(h(r_i,t_j)) = q_2(u(r_i,t_j),t_j) =
q_2(h_j(r_i),t_j) = q_{2j}(h_j(r_i))$ if $q_{2j}(r)$ is interpreted to be $q_2(r,t_j)$ for every
$r\in [0,1]$. Thus, with the trapezoidal rule the integral is discretized by
\[ E(h_j) = \frac{1}{2}\sum_{i=1}^{M-1} \Delta r_i(E_{i+1}^j+E_i^j) \]
where for $i=1,\ldots,M$,
\[ E_i^j=||q_{1j}(r_i) - q_{2j}(h_j(r_i))\sqrt{h_j'(r_i)}||^2. \]
From this, again using the trapezoidal rule, we can then discretize the double integral
 \[ \int_0^1\int_0^1 ||q_1(r,t) -\big(\; (q_2\circ h)\sqrt{\frac{\partial (u,v)}{\partial (r,t)}}\;\big)
(r,t)||^2dr\,dt \]
by
\[ E = \frac{1}{2}\sum_{j=1}^{N-1} \Delta t_j(E(h_{j+1})+E(h_j)), \]
where for $j=1,\ldots,N-1$, $\Delta t_j = t_{j+1}-t_j$.
\\ \smallskip\\
Given $j$, $1\leq j\leq N$, treating now $h_j(r_i)$, $i=1,\ldots,M$, in the definition of $E(h_j)$ as the
discretization of any diffeomorphism $h_j$ from $[0,1]$ onto $[0,1]$, if for each $j$, $j=1,\ldots,N$,
we can find $h_j$ whose discretization minimizes $E(h_j)$, then the collection of diffeomorphisms $h_j$,
$j=1,\ldots,N$, minimizes~$E$, and a homeomorphism $h$ in $\Sigma_1$ can be identified such that
$h(r_i,t_j)=(h_j(r_i),t_j)$, $i=1,\ldots,M$, $j=1,\ldots,N$. Thus, the double integral above is
approximately minimized by $h$ among all homeomorphisms in~$\Sigma_1$, with the value of the double
integral approximately equal to $E$.
\\ \smallskip\\
In \cite{bernal4}, algorithm {\em adapt-DP}, an algorithm based on Dynamic Programming, was presented for 
approximately computing, ignoring rotations, the elastic shape registration of two curves in $d-$dimensional
space. The algorithm was originally presented in \cite{bernal} for $d=2$. Given that $\hat{q}_1$ and
$\hat{q}_2$ are discretizations of the shape functions of the two curves, $\hat{q}_1$ and
$\hat{q}_2$ are used as input for algorithm {\em adapt-DP} to compute a discretization of a diffeomorphism
for reparametrizing the second curve, the reparametrization then resulting in an approximate elastic shape
registration of the two curves. Even though for $j=1,\ldots,N$, $q_{1j}$ and $q_{2j}$ as defined above are not
exactly computed as discretizations of the shape functions of curves in $3-$dimensional space, with algorithm
{\em adapt-DP} for $d=3$ with $q_{1j}$, $q_{2j}$ taking the place of $\hat{q}_1$, $\hat{q}_2$, respectively,
we can still compute the discretization of some diffeomorphism $h_j$, i.e., $h_j(r_i)$, $i=1,\ldots,M$, that
approximately minimizes $E(h_j)$. Having done this for each $j$, $j=1,\ldots,N$, $h$ in $\Sigma_1$ can then
be identified such that $h(r_i,t_j)=(h_j(r_i),t_j)$, $i=1,\ldots,M$, $j=1,\ldots,N$, and, ignoring rotations,
$c_1$, $c_2(h)$ are interpreted to achieve approximately the partial elastic shape registration of the two
surfaces. Computing $E = \frac{1}{2}\sum_{j=1}^{N-1} \Delta t_j(E(h_{j+1})+E(h_j))$, again ignoring rotations,
then $\sqrt{E}$ is interpreted to be approximately the elastic shape distance between the two surfaces
corresponding to the partial elastic shape registration of the two surfaces.
\section{\large Computation of Rotation Matrix for Rigid Alignment of Surfaces}
In this section, we describe the computation of an approximately optimal rotation matrix for the rigid
alignment of two simple surfaces of unit area in $3-$dimensional space.
Given that $S_1$ and $S_2$ are the two surfaces, with $D$, $c_1$, $c_2$, $q_1$, $q_2$ as in the previous
section, we hope to minimize
\[ \int\int_D ||q_1(r,t)-Rq_2(r,t)||^2dr\,dt = \int_0^1\int_0^1 ||q_1(r,t)-Rq_2(r,t)||^2dr\,dt \]
with respect to rotation matrices $R$ in $3-$dimensional space, i.e., with respect to $3\times 3$ matrices
$R$ that are orthogonal and have determinant equal to~1, i.e., with respect to matrices $R$ in~$SO(3)$.
\\ \smallskip \\
As in the previous section, we must work with a discretized version of the problem. Thus we assume again that
for positive integers $M$, $N$, not necessarily equal, and partitions of $[0,1]$, $\{r_i\}_{i=1}^M$,
$r_1=0<r_2<\ldots <r_M=1$, $\{t_j\}_{j=1}^N$, $t_1=0<t_2<\ldots <t_N=1$, not necessarily uniform, $c_1$ and
$c_2$ are given as lists of $M\times N$ points in the surfaces $S_1$ and $S_2$, respectively, the lists for
$c_1$ and $c_2$ corresponding to $c_1(r_i,t_j)$ and $c_2(r_i,t_j)$, $i=1,\ldots,M$, $j=1,\ldots,N$,
respectively, and that $q_1(r_i,t_j)$ and $q_2(r_i,t_j)$ are approximately computed from $c_1(r_i,t_j)$ and
$c_2(r_i,t_j)$, $i=1,\ldots,M$, $j=1,\ldots,N$, as in the previous section.
With $\Delta r_i = r_{i+1}-r_i$, $i=1,\ldots,M-1$, for $R$ in $SO(3)$, and an integer $j$, $1\leq j \leq N$,
next with the trapezoidal rule we discretize the integral
\[ \int_0^1 ||q_1(r,t_j)-Rq_2(r,t_j)||^2dr \]
by
\begin{eqnarray*}
 F_j &=& \frac{1}{2} \sum_{i=1}^{M-1}\Delta r_i(\|q_1(r_i,t_j)-Rq_2(r_i,t_j)\|^2
+ \|q_1(r_{i+1},t_j)-Rq_2(r_{i+1},t_j)\|^2)\\
\ &=& \sum_{i=1}^M\Delta\tilde{r}_i\,\| q_1(r_i,t_j) - R q_2(r_i,t_j) \|^2,
\end{eqnarray*}
where $\Delta\tilde{r}_1=(r_2-r_1)/2$, $\Delta\tilde{r}_M=(r_M-r_{M-1})/2$, and for \mbox{$i=2,\ldots,M-1$},
$\Delta\tilde{r}_i=(r_{i+1}-r_{i-1})/2$. Note, $\Delta\tilde{r}_i>0$ for $i=1,\ldots,M$, and
$\sum_{i=1}^M\Delta\tilde{r}_i=1$.\\ \smallskip\\
From this, with $\Delta t_j = t_{j+1}-t_j$, $j=1,\ldots,N-1$, again using the trapezoidal rule and noting
that $\| Rq_2(r_i,t_j)\| = \| q_2(r_i,t_j)\|$, $i=1,\ldots,M$, $j=1,\ldots,N$, we can then discretize the
double integral
\[ \int_0^1\int_0^1 ||q_1(r,t)-Rq_2(r,t)||^2dr\,dt \]
by
\begin{eqnarray*}
F &=& \frac{1}{2}\sum_{j=1}^{N-1} \Delta t_j(F_{j+1}+F_j) = \sum_{j=1}^N \Delta\tilde{t}_j F_j\\
\ &=& \sum_{j=1}^N \Delta\tilde{t}_j
(\sum_{i=1}^M\Delta\tilde{r}_i\,\| q_1(r_i,t_j) - R q_2(r_i,t_j) \|^2)\\
\ &=& \sum_{j=1}^N \Delta\tilde{t}_j
(\sum_{i=1}^M\Delta\tilde{r}_i(\| q_1(r_i,t_j) \|^2 + \| q_2(r_i,t_j) \|^2))\\
\ &\ & -2\sum_{j=1}^N \Delta\tilde{t}_j
(\sum_{i=1}^M\Delta\tilde{r}_i((q_1(r_i,t_j))^T R q_2(r_i,t_j))),
\end{eqnarray*}
where $\Delta\tilde{t}_1=(t_2-t_1)/2$, $\Delta\tilde{t}_N=(t_N-t_{N-1})/2$, and for \mbox{$j=2,\ldots,N-1$},
$\Delta\tilde{t}_j=(t_{j+1}-t_{j-1})/2$. Note, $\Delta\tilde{t}_j>0$ for $j=1,\ldots,N$, and
$\sum_{j=1}^N\Delta\tilde{t}_j=1$.\\ \smallskip\\
Thus, minimizing $F$ over all rotations $R$ in $SO(3)$ is equivalent to maximing
over the same set of rotations
\[ \sum_{j=1}^N \Delta\tilde{t}_j
(\sum_{i=1}^M\Delta\tilde{r}_i((q_1(r_i,t_j))^T R q_2(r_i,t_j)))=\mathrm{tr}(RA^T), \]
where $A$ is the $3\times 3$ matrix with entries
\[ A_{kl}= \sum_{j=1}^N \Delta\tilde{t}_j (\sum_{i=1}^M\Delta\tilde{r}_i(q_1(r_i,t_j)_k  q_2(r_i,t_j)_l)), \]
for each pair $k,l = 1,2,3$, 
$q_1(r_i,t_j)_k$ the $k^{th}$ coordinate of $q_1(r_i,t_j)$, and $q_2(r_i,t_j)_l$ the $l^{th}$ coordinate of
$q_2(r_i,t_j)$, $i=1,\ldots,N$, $j=1,\ldots,M$, and $\mathrm{tr}(RA^T)$ is the trace of the matrix~$RA^T$.
\\ \smallskip\\
Accordingly, an optimal rotation matrix $R$ for maximizing $\mathrm{tr}(RA^T$) can be
computed from the singular value decomposition of $A$ or, more precisely, with the
Kabsch-Umeyama algorithm \cite{kabsch1,kabsch2,umeyama,lawrence,bernal3}
(see Algorithm Kabsch-Umeyama below for $3-$dimensional surfaces, where
$\mathrm{diag}\{s_1,s_2,s_3\}$
is the $3\times 3$ diagonal matrix with numbers $s_1,s_2,s_3$ as the elements
of the diagonal, in that order running from the upper left to the lower right of
the matrix). A {\em singular value decomposition} (SVD)~\cite{lay} of~$A$ is a
representation of the form $A=USV^T$, where $U$ and $V$ are~$3\times 3$ orthogonal
matrices and $S$ is a $3\times 3$ diagonal matrix with the singular values of $A$, which
are nonnegative real numbers, appearing in the diagonal of $S$ in descending order,
from the upper left to the lower right of~$S$. Finally, note that the SVD concept can be
generalized so that any matrix of any dimension, not
necessarily square, has a singular value decomposition, not necessarily unique~\cite{lay}.
\begin{algorithmic}
\STATE \noindent\rule{13cm}{0.4pt}
\STATE {\bf Algorithm Kabsch-Umeyama} for surfaces (KU3 algorithm)
\STATE \noindent\rule[.1in]{13cm}{0.4pt}
\STATE Set $\Delta\tilde{r}_1=(r_2-r_1)/2$, $\Delta\tilde{r}_M=(r_M-r_{M-1})/2$, and for \mbox{$i=2,\ldots,M-1$},
$\Delta\tilde{r}_i=(r_{i+1}-r_{i-1})/2$.
\STATE Set $\Delta\tilde{t}_1=(t_2-t_1)/2$, $\Delta\tilde{t}_N=(t_N-t_{N-1})/2$, and for \mbox{$j=2,\ldots,N-1$},
$\Delta\tilde{t}_j=(t_{j+1}-t_{j-1})/2$.
\STATE Set $q_1(r_i,t_j)_k$ equal to the $k^{th}$ coordinate of $q_1(r_i,t_j)$ for $i=1,\ldots,M$,
$j=1,\ldots,N$, $k=1,2,3$.
\STATE Set $q_2(r_i,t_j)_l$ equal to the $l^{th}$ coordinate of $q_2(r_i,t_j)$ for $i=1,\ldots,M$,
$j=1,\ldots,N$, $l=1,2,3$.
\STATE Compute $A_{kl}= \sum_{j=1}^N \Delta\tilde{t}_j (\sum_{i=1}^M\Delta\tilde{r}_i(q_1(r_i,t_j)_k
q_2(r_i,t_j)_l))$\\ for each pair $k,l = 1,2,3$. 
\STATE Identify $3\times 3$ matrix $A$ with entries $A_{kl}$ for each pair $k,l=1,2,3$.
\STATE Compute SVD of $A$, i.e., identify $3\times 3$ matrices $U$, $S$, $V$, so that
\STATE $A = U S V^T$ in the SVD sense.
\STATE Set $s_1=s_2=1$.
\STATE {\small\bf if} $\det(UV) > 0$ {\small\bf then} set $s_3=1$.
\STATE {\small\bf else} set $s_3=-1$. {\small\bf end if}
\STATE Set $\tilde{S} = \mathrm{diag}\{s_1,s_2,s_3\}$.
\STATE Compute and return $R = U \tilde{S} V^T$ and $maxtrace=\mathrm{tr}(RA^T)$.
\STATE \noindent\rule{13cm}{0.4pt}
\end{algorithmic}
\section{\large Procedure for Optimizing over both Rotations and Reparametrizations using Dynamic Programming}
With $D=[0,1]\times[0,1]$, $c_1$, $c_2$, $q_1$, $q_2$, $S_1$, $S_2$ as above, $R$ in $SO(3)$, $h$ in $\Sigma_0$,
so that $\frac{\partial (u,v)}{\partial (r,t)}$, the determinant of the Jacobian of~$h$, is positive on~$D$,
we hope to minimize
 \[ \int_0^1\int_0^1 ||q_1(r,t) -\big(\;R(q_2\circ h)\sqrt{\frac{\partial (u,v)}{\partial (r,t)}}\;\big)
(r,t)||^2dr\,dt \]
with respect to $R$ and $h$.\\
We note, using arguments as those in \cite{dogan}, the above minimization problem
can be reformulated as that of minimizing
 \[ \int_0^1\int_0^1 ||Rq_1(r,t) -\big(\;(q_2\circ h)\sqrt{\frac{\partial (u,v)}{\partial (r,t)}}\;\big)
(r,t)||^2dr\,dt \]
with respect to $R$ and $h$.\\
This allows for the second surface to be reparametrized while the first one is rotated. Of course, as already
noted above, we work with $\Sigma_1$, as defined above, instead of $\Sigma_0$ of which it is a subset, as this
allows for the use of Dynamic Programming when optimizing over reparametrizations of the second surface.
Assuming $M$, $N$, $r_i$, $i=1,\ldots,M$, $t_j$, $j=1,\ldots,N$, $c_1(r_i,t_j)$, $c_2(r_i,t_j)$, $q_1(r_i,t_j)$,
$q_2(r_i,t_j)$, $i=1,\ldots,M$, $j=1,\ldots,N$, are as in the previous sections, for the purpose of
approximately minimizing the second double integral above with respect to $R$ in $SO(3)$, $h$ in $\Sigma_1$,
we use the procedure below that alternates computations of discretizations of approximately optimal
homeomorphisms in $\Sigma_1$ using Dynamic Programming (one per iteration
for reparametrizing the second surface) and approximately optimal rotation matrices (one per iteration
for rotating the first surface), these computations as described in the previous two sections.
The procedure, Procedure DP-surface-min, with
$c_1(r_i,t_j)$, $c_2(r_i,t_j)$, $q_1(r_i,t_j)$, $q_2(r_i,t_j)$, $i=1,\ldots,M$, $j=1,\ldots,N$,
as input, is summarized below. In it, given discretizations $q(r_i)$, $\hat{q}(r_i)$, $i=1,\ldots,M$,
of functions $q$, $\hat{q}$, treated as discretizations of shape functions of curves in $3-$dimensional space,
to say ``Execute DP algorithm for $q(r_i)$, $\hat{q}(r_i)$, $i=1,\ldots,M$'' will mean the DP algorithm
({\em adapt-DP} for~$d=3$) should be executed with $q(r_i)$, $\hat{q}(r_i)$, $i=1,\ldots,M$, as input, as
described in Section~5 above. Also, given $\tilde{q}_1(r_i,t_j)$, $\tilde{q}_2(r_i,t_j)$,
$i=1,\ldots,M$, $j=1,\ldots,N$, discretizations of shape functions $\tilde{q}_1$, $\tilde{q}_2$ of the two surfaces,
to say ``Execute KU3 algorithm for $\tilde{q}_1(r_i,t_j)$, $\tilde{q}_2(r_i,t_j)$, $i=1,\ldots,M$, $j=1,\ldots,N$''
will mean the Kabsch-Umeyama algorithm for surfaces, outlined in the previous section, should be executed with
$\tilde{q}_1$, $\tilde{q}_2$ taking the place of $q_1$, $q_2$, respectively, in the~algorithm.
%
\begin{algorithmic}
\STATE \noindent\rule{13cm}{0.4pt}
\STATE {\bf Procedure DP-surface-min}
\STATE \noindent\rule[.1in]{13cm}{0.4pt}
\STATE Set $\Delta r_i = r_{i+1}-r_i$, $i=1,\ldots,M-1$.
\STATE Set $\Delta\tilde{r}_1=(r_2-r_1)/2$, $\Delta\tilde{r}_M=(r_M-r_{M-1})/2$, and for \mbox{$i=2,\ldots,M-1$},
$\Delta\tilde{r}_i=(r_{i+1}-r_{i-1})/2$.
\STATE Set $\Delta\tilde{t}_1=(t_2-t_1)/2$, $\Delta\tilde{t}_N=(t_N-t_{N-1})/2$, and for \mbox{$j=2,\ldots,N-1$},
$\Delta\tilde{t}_j=(t_{j+1}-t_{j-1})/2$.
\STATE Set $\hat{q}_2(r_i,t_j)=q_2(r_i,t_j)$ for $i=1,\ldots,M$, $j=1,\ldots,N$.
\STATE Set $iter=0$, $E^{curr}=10\,^6$, $iten=10$, $tol=10\,^{-6}$.
\REPEAT
\STATE Set $iter=iter+1$, $E^{prev} = E^{curr}$.
\STATE Execute KU3 algorithm for $\hat{q}_2(r_i,t_j)$, $q_1(r_i,t_j)$, $i=1,\ldots,M$,\\ $j=1,\ldots,N$, to get
rotation matrix $R$.
\STATE Set $\hat{q}_1(r_i,t_j)=Rq_1(r_i,t_j)$ for $i=1,\ldots,M$, $j=1,\ldots,N$.
\FOR {$j=1,\ldots,N$}
\STATE Set $q_{1j}(r_i)=\hat{q}_1(r_i,t_j)$, $q_{2j}(r_i)=q_2(r_i,t_j)$, $i=1,\ldots,M$.
\STATE Execute DP algorithm for $q_{1j}(r_i)$, $q_{2j}(r_i)$, $i=1,\ldots,M$, to get
\STATE discretization of diffeomorphism $h_j$:  $h_j(r_i)$, $i=1,\ldots,M$.
\STATE Set $h_j'(r_i) = (h_j(r_{i+1})-h_j(r_i))/\Delta r_i$ for $i=1,\ldots,M-1$,\\
       $h_j'(r_M) = h_j'(r_1)$.
\STATE From interpolation of $q_{2j}(r_i)$, $i=1,\ldots,M$, with a cubic spline\\
       set $\hat{q}_{2j}(r_i)= \sqrt{h_j'(r_i)}q_{2j}(h_j(r_i))$ for $i=1,\ldots,M$.
\STATE Compute $E(h_j)=
\sum_{i=1}^M\Delta\tilde{r}_i\|q_{1j}(r_i)-\sqrt{h_j'(r_i)}q_{2j}(h_j(r_i))\|^2$
\STATE $=\sum_{i=1}^M\Delta\tilde{r}_i\|q_{1j}(r_i)-\hat{q}_{2j}(r_i)\|^2$.
\STATE Set $\hat{q}_2(r_i,t_j)=\hat{q}_{2j}(r_i)$ for $i=1,\ldots,M$.
\ENDFOR
\STATE $E^{curr} =\sum_{j=1}^N\Delta\tilde{t}_jE(h_j)$.
\UNTIL {$|E^{curr}-E^{prev}|< tol$ or $iter > iten$.}
\STATE $E=E^{curr}$.
\FOR {$j=1,\ldots,N$}
\STATE Set $h(r_i,t_j)=(h_j(r_i),t_j)$ for $i=1,\ldots,M$.
\STATE Set $c_{2j}(r_i)=c_2(r_i,t_j)$ for $i=1,\ldots,M$.
\STATE From interpolation of $c_{2j}(r_i)$, $i=1,\ldots,M$, with a cubic spline\\
       set $\hat{c}_2(r_i,t_j)= c_{2j}(h_j(r_i))$ for $i=1,\ldots,M$.
\ENDFOR
\STATE Set $\hat{c}_1(r_i,t_j)=Rc_1(r_i,t_j)$ for $i=1,\ldots,M$, $j=1,\ldots,N$.
\STATE Return $E$, $R$, $h(r_i,t_j)$,
$\hat{c}_1(r_i,t_j)$, $\hat{q}_1(r_i,t_j)$, $\hat{c}_2(r_i,t_j)$, $\hat{q}_2(r_i,t_j)$,\\
$i=1,\ldots,M$, $j=1,\ldots,N$.
\STATE \noindent\rule{13.2cm}{0.4pt}
\end{algorithmic}
On output, restricting ourselves to homeomorphisms in $\Sigma_1$, $E$ is interpreted to be the square of the
elastic shape distance between $c_1$ and $c_2$; $\hat{c}_1(r_i,t_j)$ and $\hat{c}_2(r_i,t_j)$, $i=1,\ldots,M$,
$j=1,\ldots,N$, are interpreted to achieve the elastic shape registration of $c_1$ and $c_2$; $\hat{q}_1$ and
$\hat{q}_2$ are the shape functions of $\hat{c}_1$ and $\hat{c}_2$, respectively; $R$ is the optimal rotation
matrix and $h$ is the optimal homeomorphism in~$\Sigma_1$ with which everything is computed.
Everything including $R$ and $h$ approximately computed.
\section{\large Results from Implementation of Methods}
\begin{figure}
\centering
\begin{tabular}{ccc}
\includegraphics[width=0.3\textwidth]{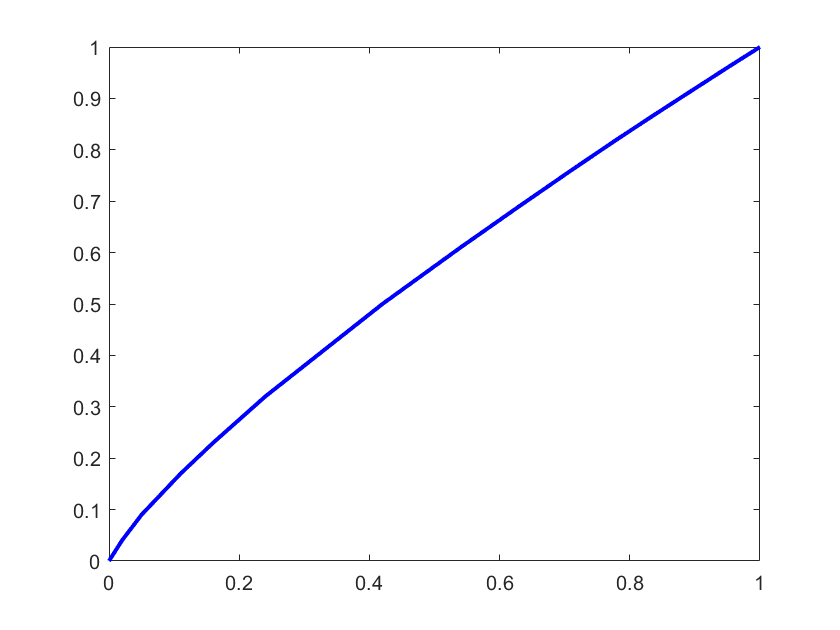}
&
\includegraphics[width=0.3\textwidth]{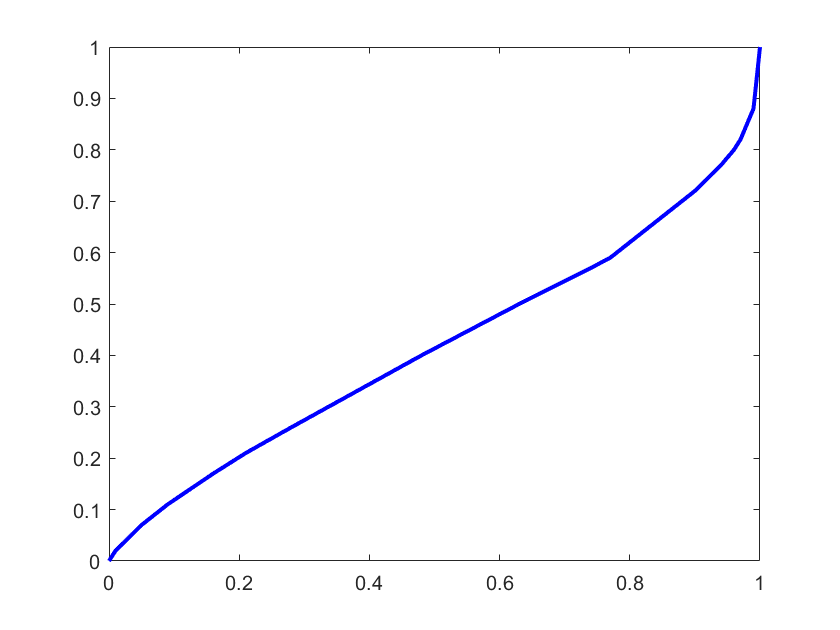}
&
\includegraphics[width=0.3\textwidth]{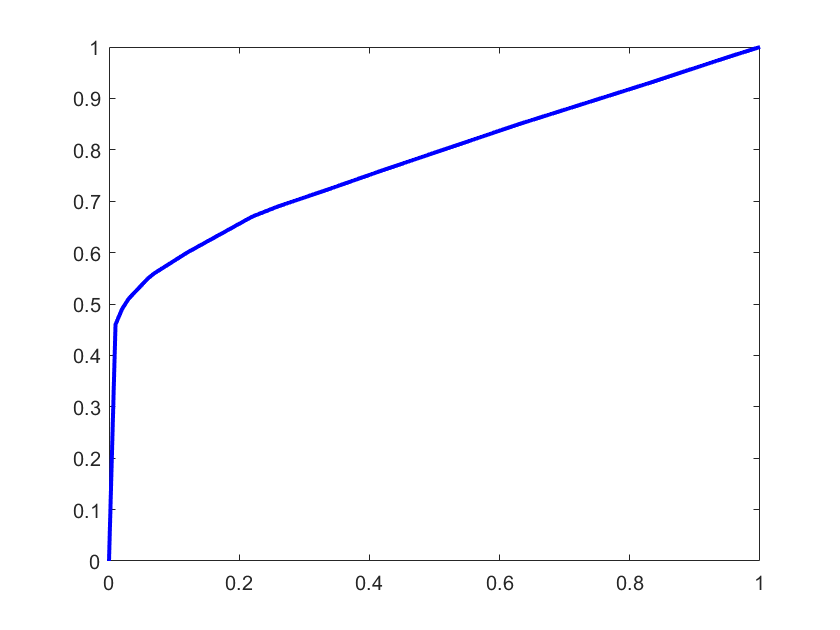}
\end{tabular}
\caption{\label{F:gammas}
Graphs of optimal diffeomorphisms from execution of Dynamic Programming software on the three pairs of
surfaces with $\gamma(r,t)=(r^{5/4},t)$, $(r,t) \in [0,1]\times [0,1]$.
One per pair, as for each pair the same two curves in $3-$d space were used as input to the
software each time it was executed. Thus the same diffeomorphism was computed each time for each pair
of surfaces.
}
\end{figure}

\begin{figure}
\centering
\begin{tabular}{ccc}
\includegraphics[width=0.3\textwidth]{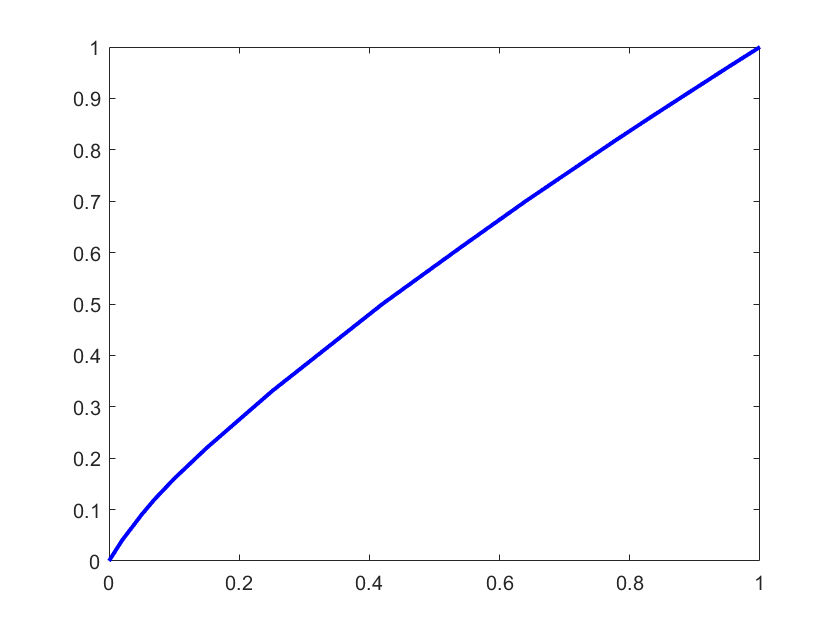}
&
\includegraphics[width=0.3\textwidth]{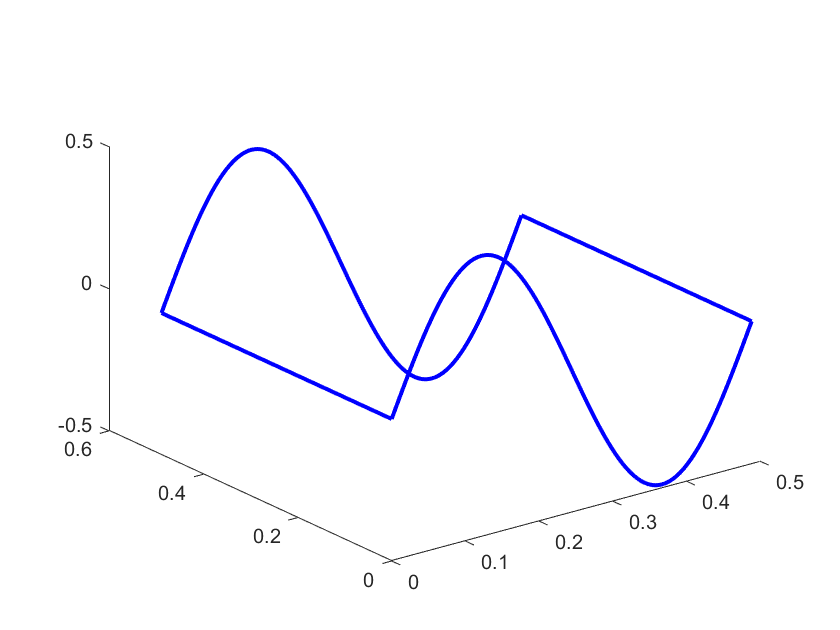}
&
\includegraphics[width=0.3\textwidth]{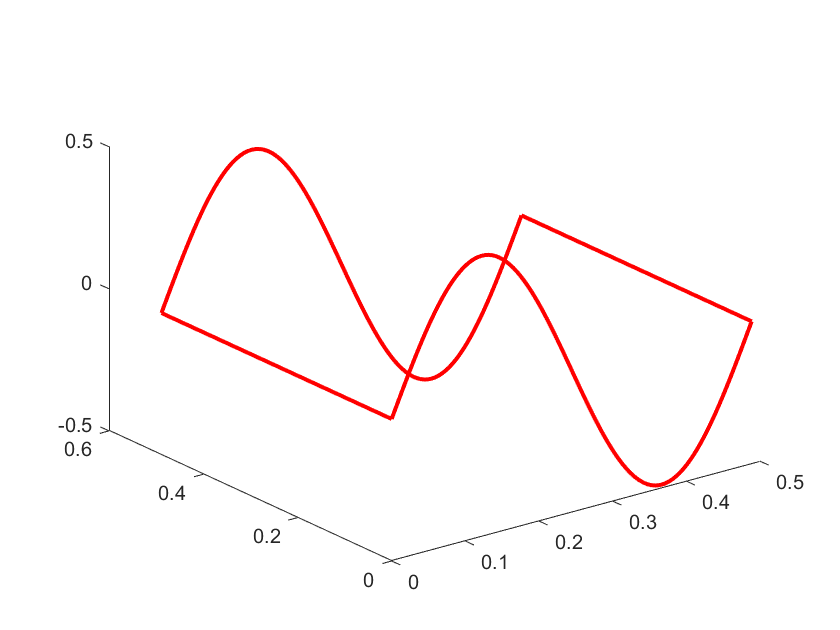}
\end{tabular}
\caption{\label{F:f121best}
For $\gamma(r,t)=(r^{5/4},t^{5/4})$, $(r,t) \in [0,1]\times [0,1]$, views of graph of optimal
diffeomorphism computed the $51^{st}$ time the Dynamic Programming software was executed
(leftmost plot), of boundary of optimally rotated first surface (middle plot), and of optimally
reparametrized second surface (rightmost plot) after computation of partial elastic shape~\mbox{registration}.
}
\end{figure}

A software package that incorporates the methods presented in this paper for computing, using Dynamic
Programming, a partial elastic shape registration of two simple surfaces in $3-$dimensional space, and
therefore the elastic shape distance between them associated with this partial registration, has been
implemented. The implementation is in Matlab\footnote{The identification of any commercial product or trade
name does not imply endorsement or recommendation by the National Institute of Standards and Technology.}
with the exception of the Dynamic Programming routine which is written in Fortran but is executed
as a Matlab mex file. In this section, we present results obtained from executions of the software
package. We note, the software package as well as input data files, a README file, etc.
can be obtained at the following link\smallskip\\
\hspace*{.35in}\verb+https://doi.org/10.18434/mds2-3056+
\smallskip\\
We note, Matlab file ESD\_\,driv\_\,surf\_\,3d.m is the driver routine of the package, and Fortran routine
DP\_\,MEX\,\_\,WNDSTRP\_\,ALLDIM.F is the Dynamic Programming routine which has already been processed
(with parameter dimx =~3) to be executed as a Matlab mex file. In case the Fortran routine must be
processed to obtain a new mex file, this can be done by typing in the Matlab window:
mex -\,compatibleArrayDims DP\_\,MEX\,\_\,WNDSTRP\_\,ALLDIM.F
\\ \smallskip\\
At the start of the execution of the software, we assume $S_1$, $S_2$ are the two simple surfaces in
$3-$dimensional space under consideration, with functions $c_1$,
$c_2: D\equiv [0,T_1]\times [0,T_2]\rightarrow \mathbb{R}^3$, $T_1$, $T_2>0$, as their parametrizations,
respectively, so that $S_1=c_1(D)$, $S_2=c_2(D)$. We also assume that as input to the software, for
positive integers $M$, $N$, not necessarily equal, and partitions of $[0,T_1]$, $[0,T_2]$, respectively,
$\{r_i\}_{i=1}^M$, $r_1=0<r_2<\ldots <r_M=T_1$, $\{t_j\}_{j=1}^N$, $t_1=0<t_2<\ldots <t_N=T_2$, not
necessarily uniform, discretizations of $c_1$, $c_2$ are given, each discretization in the form of a list
of $M\times N$ points in the corresponding surface, namely $c_1(r_i,t_j)$ and $c_2(r_i,t_j)$, $i=1,\ldots,M$,
$j=1,\ldots,N$, respectively, and for $k=1,2$, as specified in the Introduction section, in the order
$c_k(r_1,t_1)$, $c_k(r_2,t_1)$, $\ldots$,
$c_k(r_M,t_1)$, $\ldots$, $c_k(r_1,t_N)$, $c_k(r_2,t_N)$, $\ldots$, $c_k(r_M,t_N)$. Based on this input,
for the purpose of computing, using Dynamic Programming, a partial elastic shape registration of $S_1$ and
$S_2$, together with the elastic shape distance between them associated with the partial registration,
the program always proceeds first to scale the partitions $\{r_i\}_{i=1}^M$, $\{t_j\}_{j=1}^N$, so that
they become partitions of $[0,1]$, and to compute an approximation of the area of each surface.
During the execution of the software package, the former is accomplished by Matlab routine 
ESD\_\,driv\_\,surf\_\,3d.m, while the latter by Matlab routine ESD\_\,comp\_\,surf\_\,3d.m through the
computation for each $k$, $k=1,2$ of the sum of the areas of triangles with vertices $c_k(r_i,t_j)$,
$c_k(r_{i+1},t_{j+1}$, $c_k(r_i,t_{j+1})$, and $c_k(r_i,t_j)$, $c_k(r_{i+1},t_j)$, $c_k(r_{i+1},t_{j+1})$,
for $i=1,\ldots,M-1$, $j=1,\ldots,N-1$.
The program then proceeds to scale the discretizations of the parametrizations of the two
surfaces so that each surface has approximate area equal to~1 (given a surface and its approximate area, each
point in the discretization of the parametrization of the surface is divided by the square root of half the
approximate area of the surface).  Once routine ESD\_\,comp\_\,surf\_\,3d.m is done, the actual computations
of the partial registration and associated elastic shape distance are carried out by Matlab routine
ESD\_\,core\_\,surf\_\,3d.m in which the methods for this purpose presented in this paper,
mainly Procedure DP-surface-min in Section~7, have been implemented.
\\ \smallskip\\
The results that follow were obtained from applications of our software package on discretizations of three
kinds of surfaces in $3-$dimensional space that we call surfaces of the sine, helicoid and cosine-sine kind.
On input all surfaces were given as discretizations on the unit square
($[0,1]\times [0,1]$), each interval $[0,1]$ uniformly partitioned into $100$ intervals so that the unit
square was thus partitioned into $10000$ squares, each square of size $0.01\times 0.01$, their corners
making up a set of~$10201$ points.
Using the notation used at the beginning of this section, the uniform partitions of the
two $[0,1]$ intervals that define the unit square were then $\{r_i\}_{i=1}^M$, $\{t_j\}_{j=1}^N$, with
$M=N=101$, $r_{101}=t_{101}=1.0$, thus already scaled from the start as required,
and by evaluating the surfaces at the $10201$ points identified above in the order as specified above and
in the Introduction section, a discretization of each surface was obtained consisting of~$10201$ points.
Given a pair of surfaces of one of the three kinds mentioned above, and given that a partial elastic shape
registration of the two surfaces and the elastic shape distance between them associated with the partial
registration were to be computed, one surface was identified as the first surface, the other one
as the second surface (in the procedure for optimizing over rotations and reparametrizations using Dynamic
Programing as described in Section~7, Procedure DP-surface-min, the second surface is reparametrized while
the first one is rotated). For the purpose of testing the capability of the software for optimizing over
reparametrizations based on Dynamic Programming, again using the notation used at the beginning of this
section, for $\gamma$, a bijective function on the unit square to be defined below,
with $(\hat{r}_i,\hat{t}_j)=\gamma(r_i,t_j)$, $i=1,\ldots,101$, $j=1,\ldots,101$,
the second surface was reparametrized through its discretization, namely by setting $\hat{c}_2=c_2$ and
computing $c_2(r_i,t_j) =\hat{c}_2(\hat{r}_i,\hat{t}_j)$, $i=1,\ldots,101,\ j=1,\ldots,101$,
while the first surface was kept as originally defined and discretized by computing
$c_1(r_i,t_j)$, $i=1,\ldots,101,\ j=1,\ldots,101$. Given the pair of surfaces, the program then, using the
discretizations of the surfaces as just described in terms of $c_1$, $c_2$, etc., after computing
an approximation of the area of each surface and scaling each surface to have approximate area equal to~1,
proceeded to compute a partial elastic shape registration of the two surfaces and the elastic shape distance
between them associated with the partial registration.
We note that because $N$ equaled 101, during the execution the software package, the Dynamic Programming software
was executed 101 times each time the {\small\bf repeat} loop in Procedure DP-surface-min was executed.
Finally, we note again that in what follows we refer to homeomorphisms on~$[0,1]$ as diffeomorphisms in order
to distinguish them from homeomorphisms on the unit square in the plane.
\\ \smallskip\\
The first results that follow were obtained from applications of our software package on discretizations of
surfaces in $3-$dimensional space that are actually graphs of $3-$dimensional functions based on the sine curve.
Given $k$, a positive integer, one type of surface to which we refer as a surface of the sine kind (type~1) is
defined by
\[ x(r,t)= r,\ \ \ y(r,t) = t,\ \ \ z(r,t) = \sin k \pi  r,\ \ \ (r,t) \in [0,1]\times [0,1], \]
and another one (type~2) by
\[ x(r,t) = \sin k \pi  r,\ \ \ y(r,t)= r,\ \ \ z(r,t) = t,\ \ \ (r,t) \in [0,1]\times [0,1], \]
the former a rotation of the latter by applying the rotation matrix
$\left( \begin{smallmatrix} 0 & 1 & 0\\ 0 & 0 & 1\\ 1 & 0 & 0\\ \end{smallmatrix} \right)$
on the latter, thus of similar shape.
\begin{figure}
\centering
\begin{tabular}{ccc}
\includegraphics[width=0.3\textwidth]{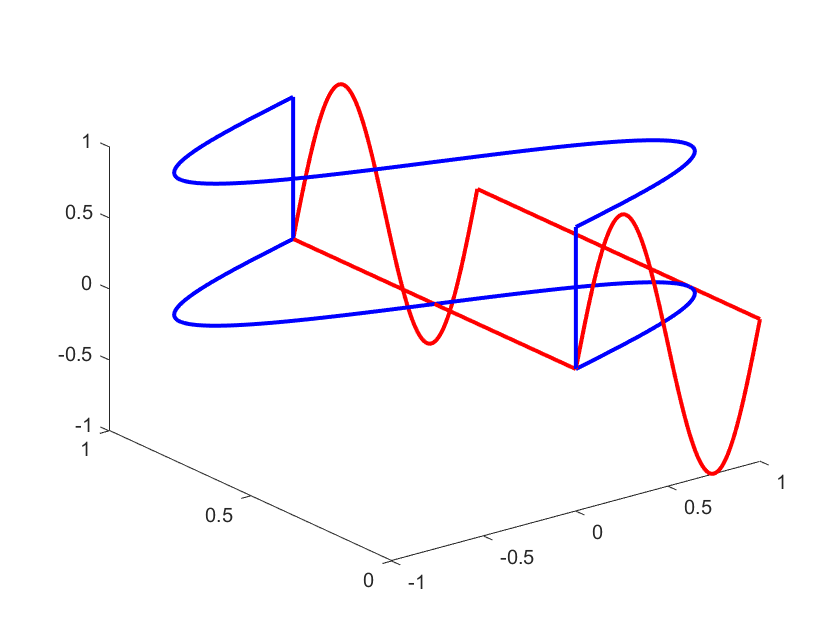}
&
\includegraphics[width=0.3\textwidth]{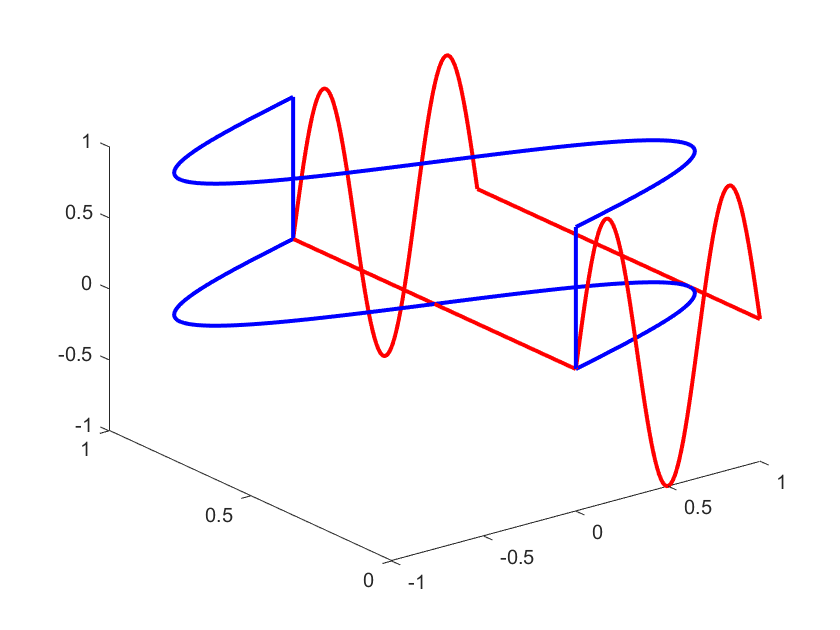}
&
\includegraphics[width=0.3\textwidth]{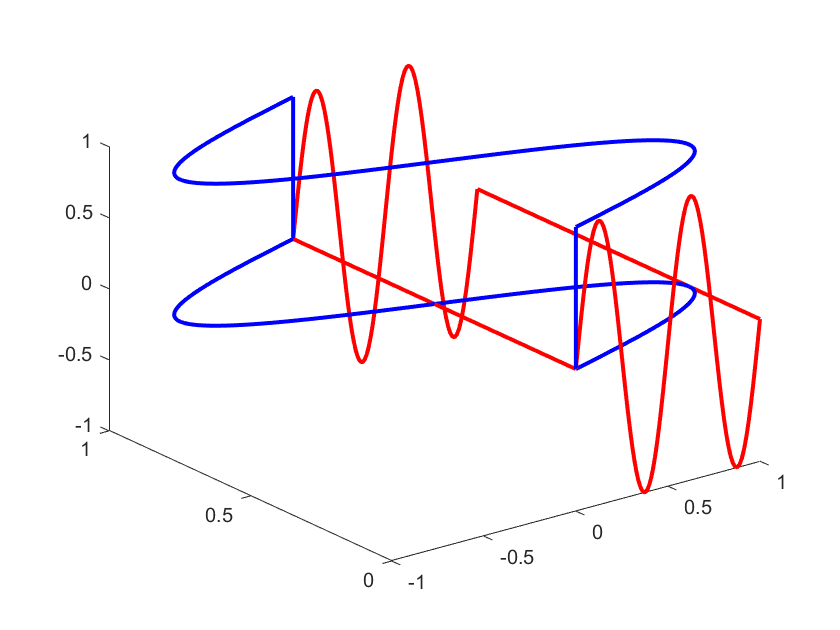}
\end{tabular}
\caption{\label{F:surfaces}
Three plots of boundaries of surfaces of the sine kind. A partial elastic shape registration of the two
surfaces in each plot and the elastic shape distance between them associated with the registration
were computed.}

\end{figure}
\\ \smallskip\\
Three plots depicting surfaces (actually their boundaries) of the sine kind for different values of $k$
are shown in Figure~\ref{F:surfaces}. (Note that in the plots there, the $x-$, $y-$, $z-$ axes are not
always to scale relative to one another). In each plot two surfaces of the sine kind appear. The two
surfaces in the leftmost plot being of similar shape, clearly the elastic shape distance between them
is exactly zero, and the hope was then that the execution of our software package applied on these two
surfaces would produce an elastic shape distance between them equal or close to zero. The type~2 surface
in each plot (in blue) was considered to be the first surface in the plot. In each plot this surface was
obtained by setting $k$ equal to~$2$ in the definition above of a type~2 surface of the sine kind so
that it is the same suface in all three plots. The other surface in each plot (in red) is a type~1
surface of the sine kind and was considered to be the second surface in each plot. From left to right in
the three plots, the second surface was obtained by setting $k$ equal to $2$, $3$, $4$, respectively, in
the definition above of a type~1 surface of the sine kind. As already mentioned above, in the procedure
for optimizing over rotations and reparametrizations using Dynamic Programing as described in Section~7,
Procedure DP-surface-min, the second surface is reparametrized while the first one is rotated.
\\ \smallskip\\
With $\gamma(r,t)=(r^{5/4},t)$, $(r,t) \in [0,1]\times [0,1]$,
all surfaces in the plots were then discretized as described above and a partial elastic shape
registration of the two surfaces in each plot and the elastic shape distance between them associated
with the partial registration were then computed through executions of our software package.
We note that for this particular $\gamma$, the discretization of the second surface was perturbed only in
the $r$ direction which made the software package more likely to succeed as Procedure DP-surface-min
always reparametrizes the second surface by applying the Dynamic Programming software exclusively on curves
in $3-$dimensional space contained in the surfaces in the $r$ direction.
The three elastic shape distances, computed in the order of the plots from left to right, were as
follows with the first distance, as hoped for, essentially equal to zero:
$\ 0.0003\ \ \ 0.3479\ \ \ 0.3192$. The times of execution in the same order
were 27, 28, 39 seconds, with the {\small\bf repeat} loop in Procedure
DP-surface-min in Section~7 executed 3, 3, 4 times, respectively. The computed optimal rotation matrix for
the pair of surfaces in the leftmost plot in Figure~\ref{F:surfaces}, was
$\left( \begin{smallmatrix} 0 & 1 & 0\\ 0 & 0 & 1\\ 1 & 0 & 0\\ \end{smallmatrix} \right).$
For the other two pairs of surfaces the computed optimal rotation matrices were both almost equal to
$\left( \begin{smallmatrix} 0 & 1 & 0\\ 0 & 0 & 1\\ 1 & 0 & 0\\ \end{smallmatrix} \right)$ as well,
their entries slightly different. It should be pointed out that because of the simplicity of surfaces of the
sine kind and the fact that for the given $\gamma$ the discretization of the second surface was perturbed
only in the $r$ direction, whenever the Dynamic Programming software was executed for a given pair of surfaces,
the same two curves in $3-$dimensional space contained in the surfaces in the $r$ direction were always used
as input to the software. Therefore for the given pair, the same solution was obtained each time (101 times)
the Dynamic Programming software was executed, in particular
the same optimal orientation-preserving diffeomorphism from $[0,1]$ onto $[0,1]$ was computed each time
together with the same elastic shape distance between the two curves in $3-$dimensional space used as input
to the software. Graphs of the optimal diffeomorphisms for each pair of surfaces in the order of the plots
from left to right in Figure~\ref{F:surfaces}, are shown in~Figure~\ref{F:gammas}.
In addition, Figure~\ref{F:f123best} shows results of the partial elastic shape registration of
the pair of surfaces in the rightmost plot in Figure~\ref{F:surfaces}.
The pair of surfaces is shown in the leftmost plot of the figure before any computations took place.
In the middle plot we see the first surface after it was rotated with the corresponding computed
optimal rotation matrix mentioned above.
In the rightmost plot we see the second surface after it was reparametrized with the homeomorphism
on the unit square corresponding to the partial elastic shape registration of the pair of surfaces, a
homeomorphism computed based on the optimal diffeomorphism obtained each time (101 times) the Dynamic
Programming software was executed for the pair of surfaces, and that because of the simplicity of the
surfaces involved and the fact that for the given $\gamma$ the discretization of the second surface
was perturbed only in the $r$ direction, was always the same diffeomorphism, the diffeomorphism whose
graph appears in the rightmost plot in~Figure~\ref{F:gammas}.
\begin{figure}
\centering
\begin{tabular}{ccc}
\includegraphics[width=0.3\textwidth]{surfaces3.png}
&
\includegraphics[width=0.3\textwidth]{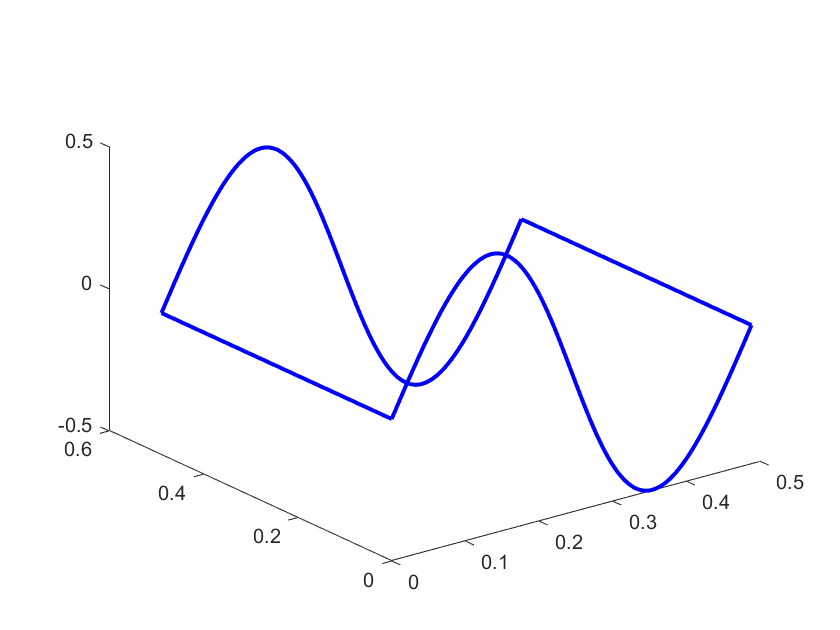}
&
\includegraphics[width=0.3\textwidth]{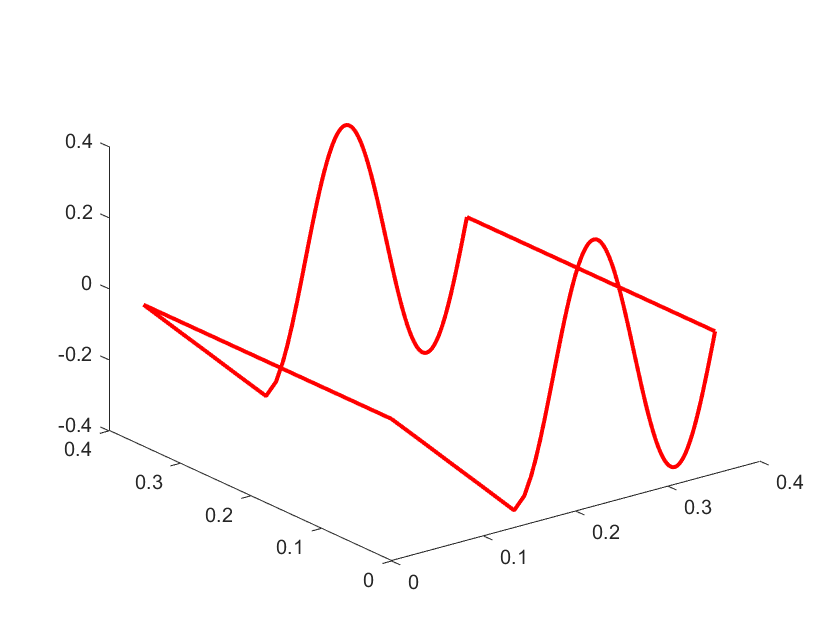}
\end{tabular}
\caption{\label{F:f123best}
With $\gamma(r,t)=(r^{5/4},t)$, $(r,t) \in [0,1]\times [0,1]$, views of boundaries of pair of surfaces
in the rightmost plot in Figure~\ref{F:surfaces} before computation of partial elastic shape registration
(leftmost plot here), of boundary of optimally rotated first surface (middle plot), and of boundary of
optimally reparametrized second surface (rightmost plot) after computations.
}

\end{figure}
\\ \smallskip\\
Finally, we note that with $\gamma(r,t)=(r^{5/4},t^{5/4})$, $(r,t) \in [0,1]\times [0,1]$,
the software package was applied on the pair of surfaces in the leftmost
plot in Figure~\ref{F:surfaces}. The computed elastic shape distance between the two surfaces was $0.0126$,
the time of execution was 28 seconds, with the {\small\bf repeat} loop in Procedure DP-surface-min in
Section~7 executed 3 times, and the computed optimal rotation matrix for the pair of surfaces was
$\left( \begin{smallmatrix} 0 & 1 & 0\\ 0 & 0 & 1\\ 1 & 0 & 0\\ \end{smallmatrix} \right).$
These results were not far from those obtained with the previous $\gamma$, however for the current $\gamma$
the discretization of the second surface was perturbed in both the $r$ and $t$ directions. As Procedure
DP-surface-min is not equipped to handle perturbations in the $t$ direction, perhaps this was the reason
why the computed elastic shape distance between the two surfaces was not exactly zero as in particular the
optimal orientation-preserving diffeomorphisms from $[0,1]$ onto $[0,1]$ computed with the Dynamic
Programming software differed slightly from one another, while the computed elastic shape distances
between the curves in $3-$dimensional space used as input to the software differed from
one another as well and were not exactly close to zero.
The graph of the optimal diffeomorphism computed the $51^{st}$ time the Dynamic Programming software
was executed is shown in Figure~\ref{F:f121best} together with results of the partial elastic shape
registration of the pair of surfaces.
It should be noted here that perhaps as long as the second surfaces
we have chosen for testing the software are perturbed in the same manner in the $r$ direction, it is
likely the graphs of the optimal diffeomorphisms computed with the Dynamic Programming software
will tend to resemble one another regardless of the surfaces involved.
\\ \smallskip\\
The next results that follow were obtained from applications of our software package on discretizations
of surfaces in $3-$dimensional space of the helicoid kind.
Given $k$, a positive integer, one type of surface to which we refer as a surface of the helicoid kind
(type~1) is defined by
\[ x(r,t)= r \cos k \pi t,\ \ \ y(r,t) = r \sin k \pi t,\ \ \ z(r,t) = k \pi t,\ \ \ (r,t) \in [0,1]\times [0,1], \]
and another one (type~2) by
\[ x(r,t) = k \pi t,\ \ \ y(r,t)= r \cos k \pi t,\ \ \ z(r,t) = r \sin k \pi t,\ \ \ (r,t) \in [0,1]\times [0,1], \]
the former a rotation of the latter by applying the rotation matrix
$\left( \begin{smallmatrix} 0 & 1 & 0\\ 0 & 0 & 1\\ 1 & 0 & 0\\ \end{smallmatrix} \right)$
on the latter, thus of similar shape.
\\ \smallskip\\
A plot depicting two surfaces (actually their boundaries) of similar shape of the helicoid
kind for $k=4$ is shown in Figure~\ref{F:helicoids}. (Note that in the plot there, the $x-$, $y-$, $z-$
axes are not always to scale relative to one another). The two surfaces being of similar shape, clearly the
elastic shape distance between them is exactly zero, and the hope was once again that the execution of our
software package applied on these two surfaces would produce an elastic shape distance between them equal
or close to zero. The type~2 surface of the helicoid kind in the plot (in blue) was considered to be
the first surface in the plot. The other surface in the plot (in red)
is a type~1 surface of the helicoid kind and was considered to be the second surface in the plot.
\begin{figure}
\centering
\begin{tabular}{ccc}
\includegraphics[width=0.4\textwidth]{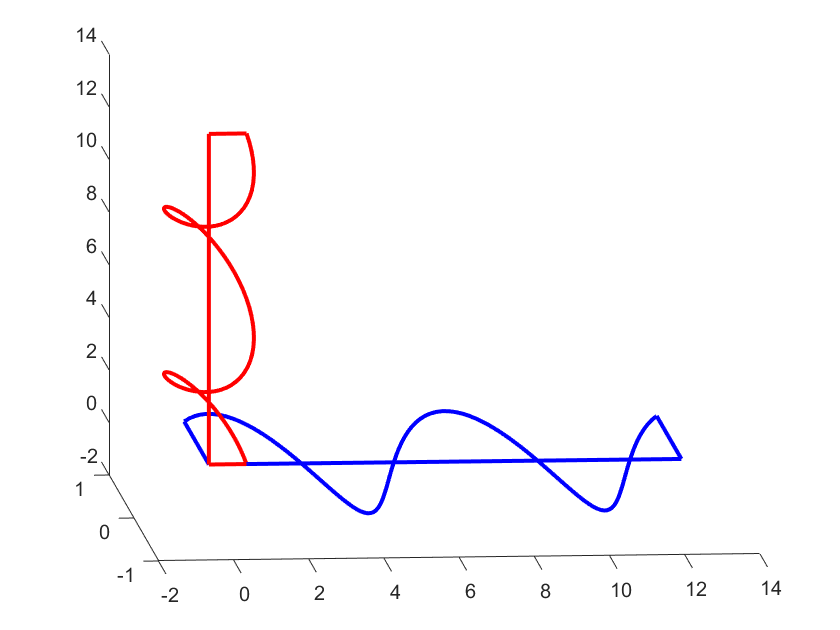}
\end{tabular}
\caption{\label{F:helicoids}
Boundaries of two surfaces of similar shape of the helicoid kind for $k=4$, type~$1$ in red, type~$2$ in blue.
}
\end{figure}
\\ \smallskip
With $\gamma(r,t)=(r^{5/4},t)$, $(r,t) \in [0,1]\times [0,1]$,
the two surfaces in the plot were then discretized as described above and a partial elastic shape
registration of the two surfaces and the elastic shape distance between them associated with the
partial registration were then computed through the execution of our software package.
Again we note that for this particular $\gamma$, the discretization of the second surface was perturbed
only in the $r$ direction which as pointed out above made the software package more likely to succeed.
The computed distance was 0.0002, which, as hoped for, was close enough to zero. The time of execution
was 15 seconds with the {\small\bf repeat} loop in Procedure DP-surface-min in Section~7 executed 2 times.
The computed optimal rotation matrix for the pair of surfaces was
$\left( \begin{smallmatrix} 0 & 1 & 0\\ 0 & 0 & 1\\ 1 & 0 & 0\\ \end{smallmatrix} \right).$
As was the case for surfaces of the sine kind, once again essentially the same solution was obtained each
time the Dynamic Programming software was executed as essentially the same two curves in $3-$dimensional
space contained in the surfaces in the $r$ direction were used each time as input to the software (the
same two curves in the sense that given a pair of curves used as input, the two curves had the same shape
and that shape was the same shape of each curve in any other pair used as input to the Dynamic Programing
software). In particular, essentially the same optimal orientation-preserving diffeomorphism from $[0,1]$
onto $[0,1]$ was computed each time together with the same elastic shape distance close to zero between the
two curves in $3-$dimensional space used as input to the software. The graph of this optimal diffeomorphism
is shown in Figure~\ref{F:h123best} together with results of the partial elastic shape registration of the
pair of~surfaces.
\begin{figure}
\centering
\begin{tabular}{ccc}
\includegraphics[width=0.3\textwidth]{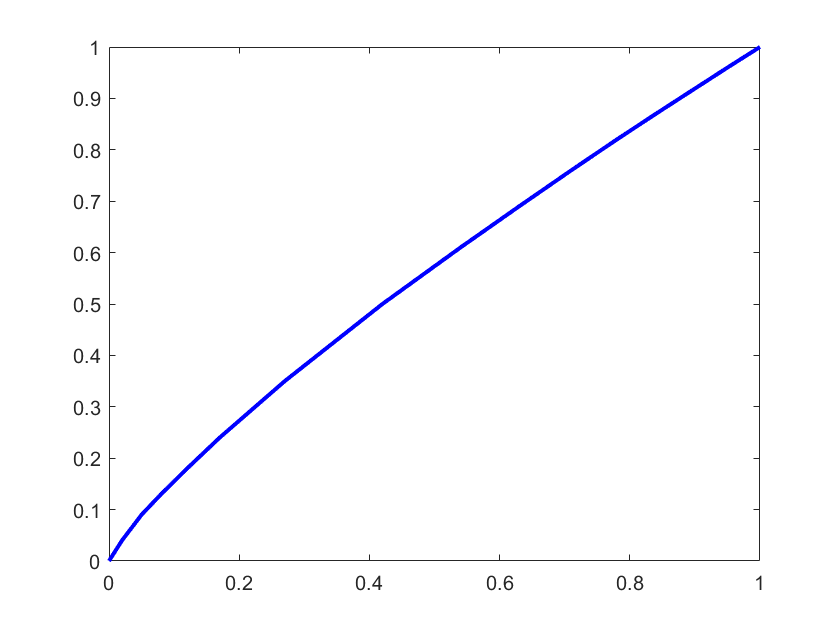}
&
\includegraphics[width=0.3\textwidth]{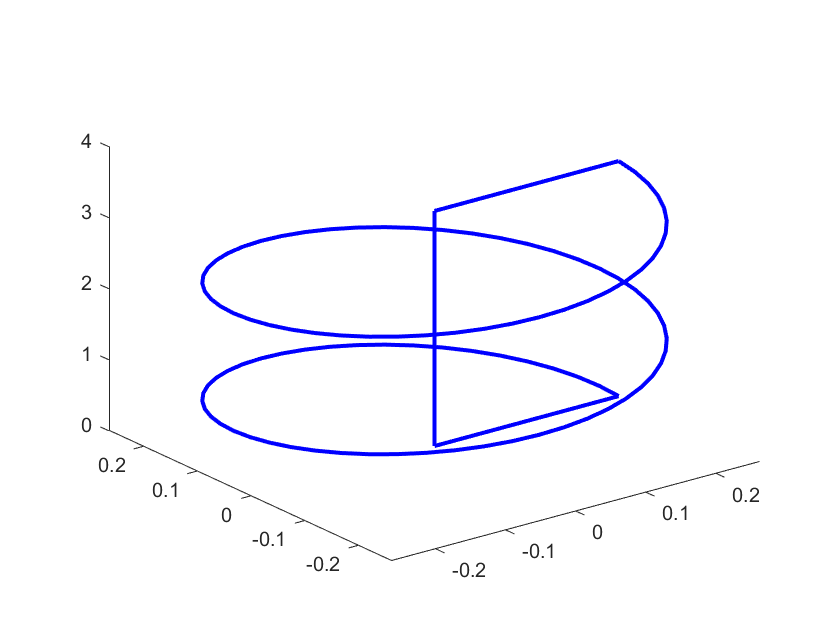}
&
\includegraphics[width=0.3\textwidth]{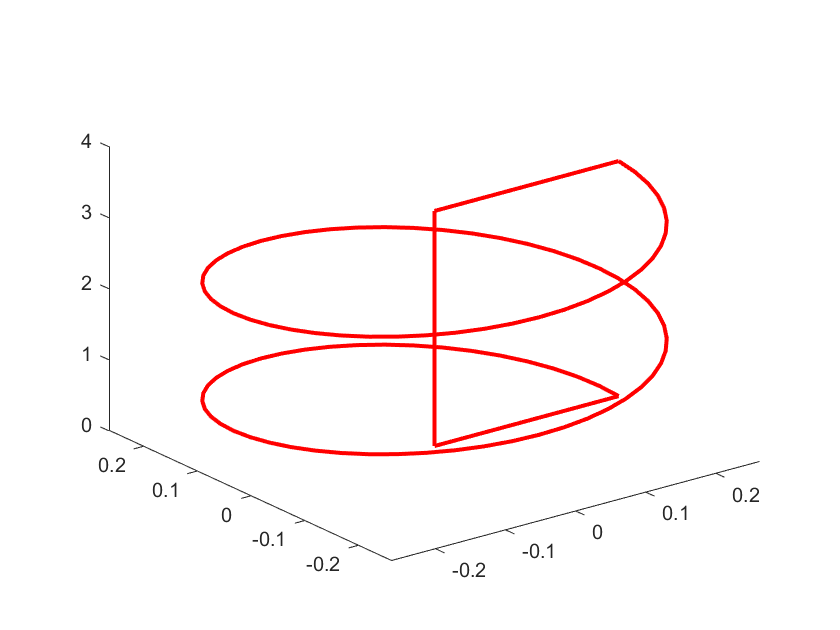}
\end{tabular}
\caption{\label{F:h123best}
For $\gamma(r,t)=(r^{5/4},t)$, $(r,t) \in [0,1]\times [0,1]$, views of graph of
optimal diffeomorphism computed each time the Dynamic Programming software was executed on pair of
surfaces (leftmost plot), of boundary of optimally rotated first surface (middle plot), and of
optimally reparametrized second surface (rightmost plot) after computation of partial elastic
shape registration.
}
\end{figure}
\\ \smallskip\\
Finally, we note that with $\gamma(r,t)=(r^{5/4},t^{5/4})$, $(r,t) \in [0,1]\times [0,1]$, the software
package was applied again on the pair of surfaces. The computed elastic shape distance between the two
surfaces was $0.0796$, the time of execution was 19 seconds, with the {\small\bf repeat} loop in
Procedure DP-surface-min in Section~7 executed 2 times, and the computed optimal rotation matrix for the
pair of surfaces was approximately
$\left( \begin{smallmatrix} .028 & .762 & .647\\ -.029 & -.646 & .763\\ .999 & -.040 & .004\\
\end{smallmatrix} \right).$
These results were not as good as those obtained with the previous $\gamma$ but still acceptable considering
that for the current $\gamma$ the discretization of the second surface was perturbed in both the $r$ and $t$
directions. As mentioned above Procedure DP-surface-min is not equipped to handle perturbations in the $t$
direction, so perhaps this was the reason why the computed elastic shape distance between the two surfaces was
not exactly zero as in particular the optimal orientation-preserving diffeomorphisms from $[0,1]$ onto
$[0,1]$ computed with the Dynamic Programming software differed slightly from one another, while the computed
elastic shape distances between the curves in $3-$dimensional space used as input to the software differed from
one another as well and were not exactly close to zero. This inability to handle perturbations in the $t$
direction may have also affected the computation of the optimal rotation matrix.
The graph of the optimal diffeomorphism computed the $51^{st}$ time the Dynamic Programming software
was executed is shown in Figure~\ref{F:h312best} together with results of the partial elastic shape
registration of the pair of surfaces.
Once again we note that perhaps as long as the second surfaces
we have chosen for testing the software are perturbed in the same manner in the $r$ direction, it is
likely the graphs of the optimal diffeomorphisms computed with the Dynamic Programming software
will tend to resemble one another regardless of the surfaces involved.
\begin{figure}
\centering
\begin{tabular}{ccc}
\includegraphics[width=0.3\textwidth]{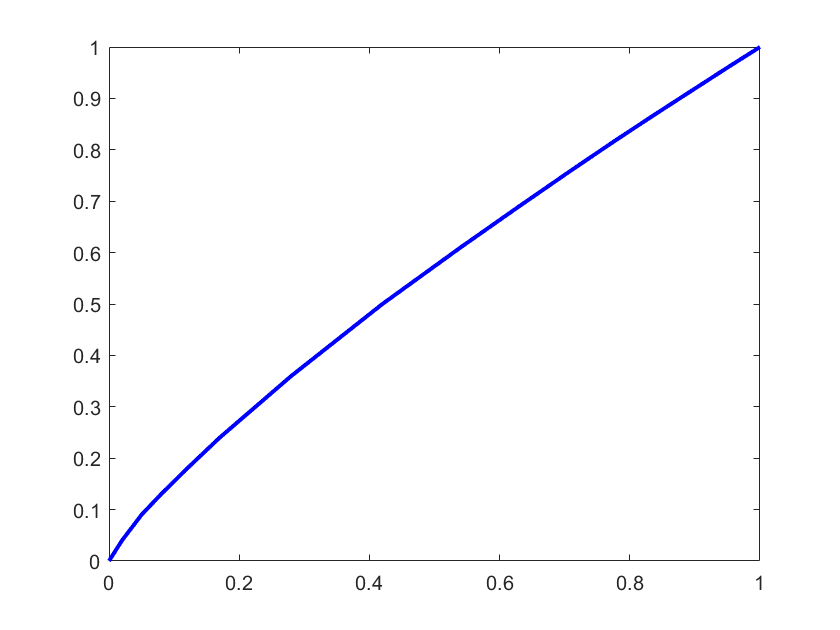}
&
\includegraphics[width=0.3\textwidth]{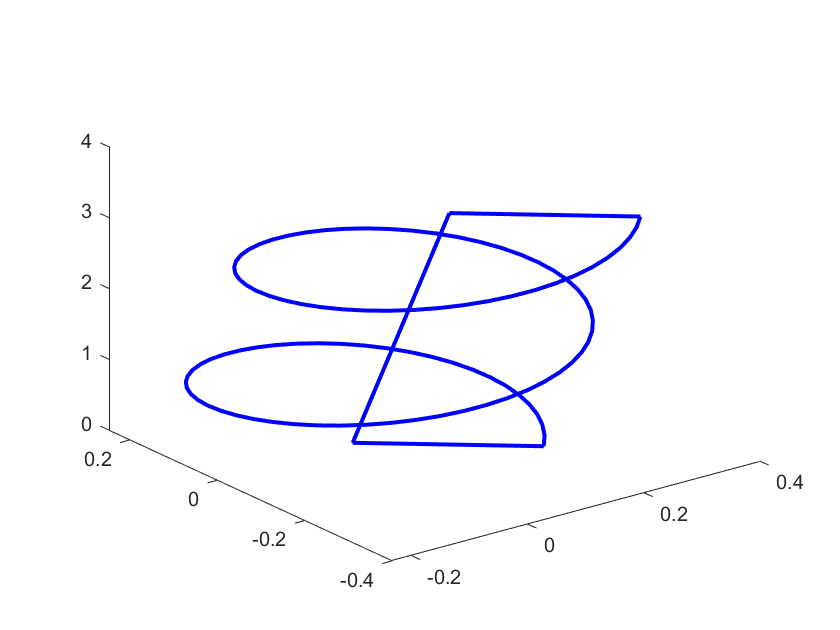}
&
\includegraphics[width=0.3\textwidth]{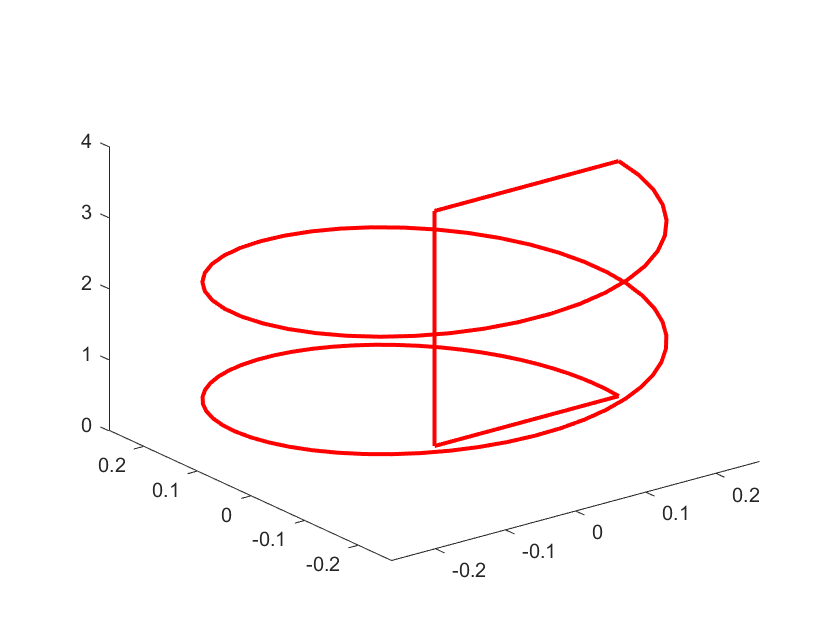}
\end{tabular}
\caption{\label{F:h312best}
For $\gamma(r,t)=(r^{5/4},t^{5/4})$, $(r,t) \in [0,1]\times [0,1]$, views of graph of
optimal diffeomorphism computed the $51^{st}$ time the Dynamic Programming software was executed
(leftmost plot), of boundary of optimally rotated first surface (middle plot), and of optimally
reparametrized second surface (rightmost plot) after computation of partial elastic shape registration.
}
\end{figure}
\\ \smallskip\\
The final results that follow were obtained from applications of our software package on discretizations of
surfaces in $3-$dimensional space that are actually graphs of $3-$dimensional functions based on the product
of the cosine and sine functions.
One surface of this kind to which we refer as a surface of the cosine-sine kind (type~1) is
defined by
\[ x(r,t)=r,\ \ \ y(r,t)=t,\ \ \ z(r,t)=(\cos0.5\pi r)(\sin0.5\pi t),\ \ \ (r,t) \in [0,1]\times [0,1], \]
and another one (type~2) by
\[ x(r,t)=(\cos0.5\pi r)(\sin0.5\pi t),\ \ \ y(r,t)=r,\ \ \ z(r,t)=t,\ \ \ (r,t) \in [0,1]\times [0,1], \]
the former a rotation of the latter by applying the rotation matrix
$\left( \begin{smallmatrix} 0 & 1 & 0\\ 0 & 0 & 1\\ 1 & 0 & 0\\ \end{smallmatrix} \right)$
on the latter, thus of similar shape.
\begin{figure}
\centering
\begin{tabular}{ccc}
\includegraphics[width=0.4\textwidth]{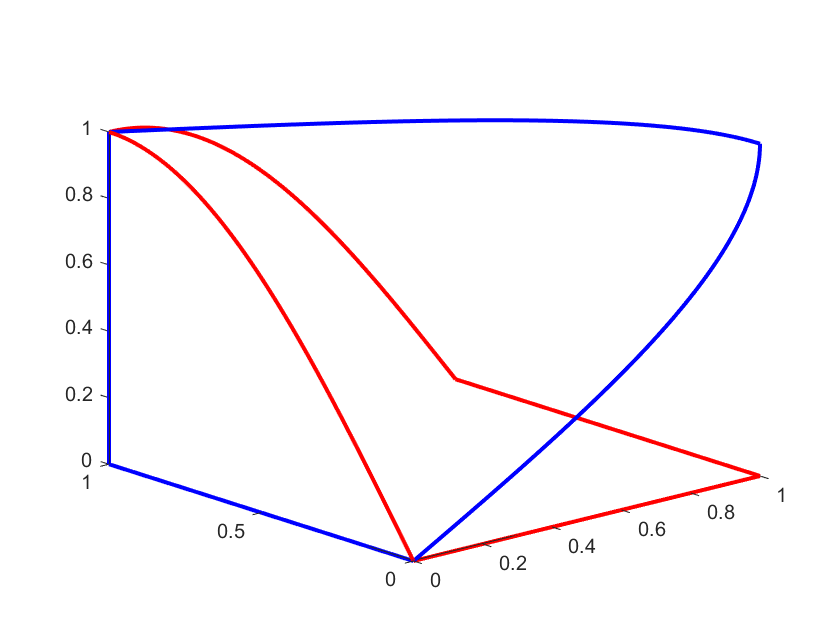}
\end{tabular}
\caption{\label{F:cos_sin}
Boundaries of two surfaces of similar shape of the cosine-sine kind, type~$1$ in red,
type~$2$ in blue.
}
\end{figure}
\\ \smallskip\\
A plot depicting two surfaces (actually their boundaries) of similar shape of the cosine-sine
kind is shown in Figure~\ref{F:cos_sin}. (Note that in the plot there, the $x-$, $y-$, $z-$
axes are not always to scale relative to one another). The two surfaces being of similar shape, clearly the
elastic shape distance between them is exactly zero, and the hope was once again that the execution of our
software package applied on these two surfaces would produce an elastic shape distance between them equal
or close to zero. The type~2 surface of the cosine-sine kind in the plot (in blue) was considered to be
the first surface in the plot. The other surface in the plot (in red)
is a type~1 surface of the cosine-sine kind and was considered to be the second surface in the plot.
\\ \smallskip\\
With $\gamma(r,t)=(r^{5/4},t)$, $(r,t) \in [0,1]\times [0,1]$,
the two surfaces in the plot were then discretized as described above and a partial elastic shape
registration of the two surfaces and the elastic shape distance between them associated with the
partial registration were then computed through the execution of our software package.
Again we note that for this particular $\gamma$, the discretization of the second surface was perturbed
only in the $r$ direction which as pointed out above made the software package more likely to succeed.
The computed distance was 0.0002, which, as hoped for, was close enough to zero. The time of execution was
22 seconds with the {\small\bf repeat} loop in Procedure DP-surface-min in Section~7 executed 3 times.
The computed optimal rotation matrix for the pair of surfaces was essentially
$\left( \begin{smallmatrix} 0 & 1 & 0\\ 0 & 0 & 1\\ 1 & 0 & 0\\ \end{smallmatrix} \right).$
It should be noted here that the type~1 surface of the cosine-sine kind satisfies that given $t_1$, $t_2$,
$0\leq t_1<t_2\leq 1$, then the two $3-$dimensional curves in the surface obtained by fixing $t$ to $t_1$
and $t$ to $t_2$, have different shapes. In spite of this, the optimal orientation-preserving diffeomorphisms
from $[0,1]$ onto $[0,1]$ computed with the Dynamic Programming software for the given $\gamma$, although
differing from one another, differed only very slightly, while the computed elastic shape distances between
the curves in $3-$dimensional space used as input to the software were all very close to zero. The graph of
the optimal diffeomorphism computed the $51^{st}$ time the Dynamic Programming software was executed is shown
in Figure~\ref{F:c123best} together with results of the partial elastic shape registration of the pair
of surfaces.
\begin{figure}
\centering
\begin{tabular}{ccc}
\includegraphics[width=0.3\textwidth]{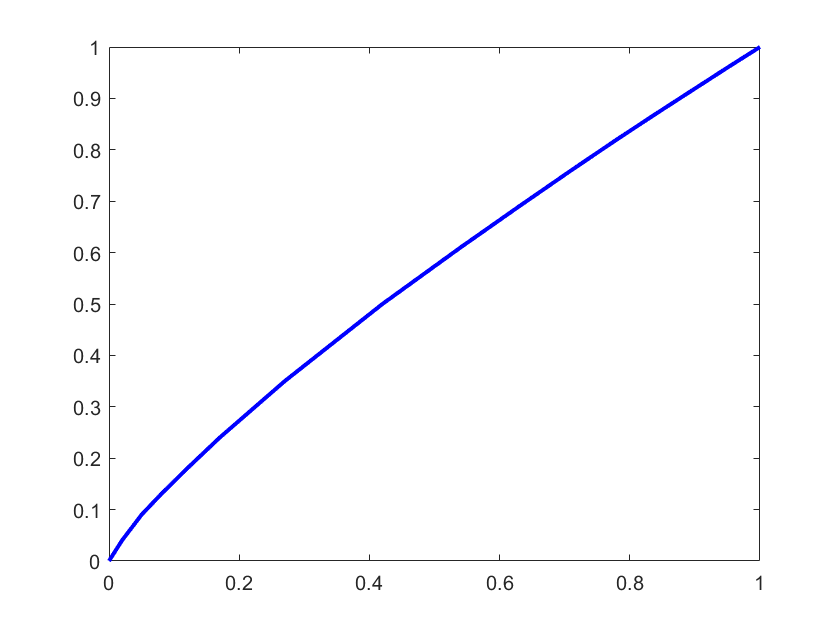}
&
\includegraphics[width=0.3\textwidth]{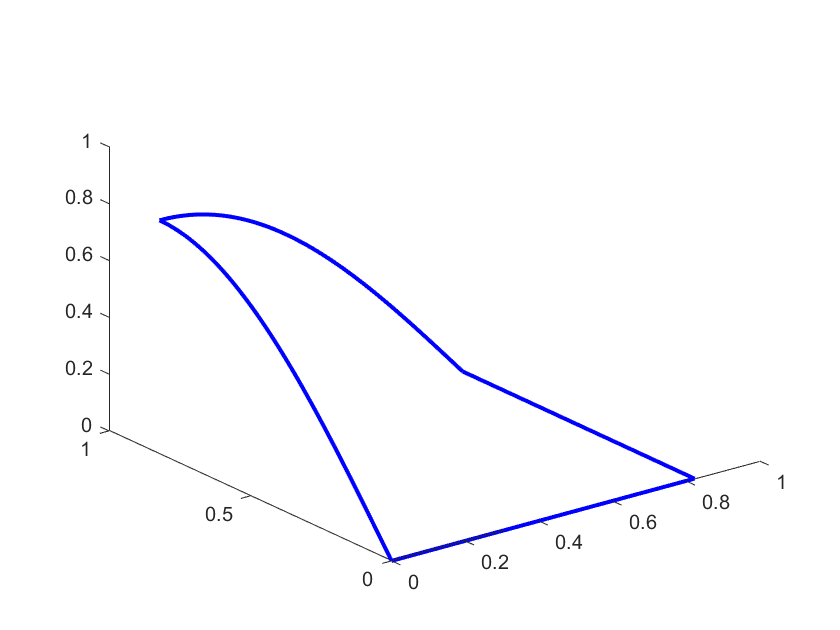}
&
\includegraphics[width=0.3\textwidth]{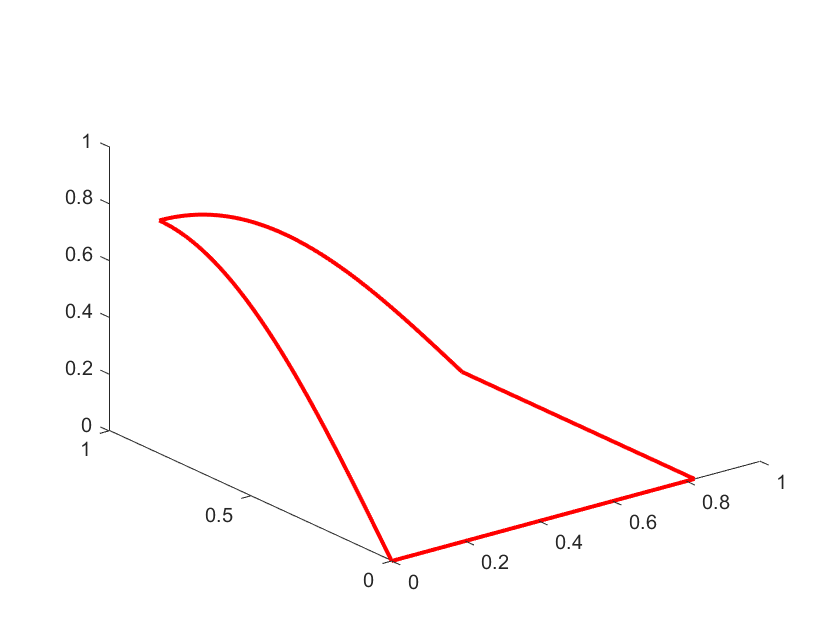}
\end{tabular}
\caption{\label{F:c123best}
For $\gamma(r,t)=(r^{5/4},t)$, $(r,t) \in [0,1]\times [0,1]$, views of graph of optimal diffeomorphism
computed the $51^{st}$ time the Dynamic Programming software was executed on pair of surfaces (leftmost
plot), of boundary of optimally rotated first surface (middle plot), and of optimally reparametrized
second surface (rightmost plot) after computation of partial registration.
}
\end{figure}
\\ \smallskip\\
Finally, we note that with $\gamma(r,t)=(r^{5/4},t^{5/4})$, $(r,t) \in [0,1]\times [0,1]$, the software
package was applied again on the pair of surfaces. The computed elastic shape distance between the two
surfaces was $0.0143$, the time of execution was 23 seconds, with the {\small\bf repeat} loop in
Procedure DP-surface-min in Section~7 executed 3 times, and the computed optimal rotation matrix for the
pair of surfaces was approximately
$\left( \begin{smallmatrix} -.043 & .999 & .026\\ -.035 & -.028 & .999\\ .998 & .042 & .036\\
\end{smallmatrix} \right).$
These results although not as good as those obtained with the previous $\gamma$ were still acceptable considering
once again that for the current $\gamma$ the discretization of the second surface was perturbed in both the $r$
and $t$ directions. Again as mentioned above Procedure DP-surface-min is not equipped to handle perturbations in the
$t$ direction, so perhaps this was the reason why the computed elastic shape distance between the two surfaces was
not exactly zero as in particular the optimal orientation-preserving diffeomorphisms from $[0,1]$ onto
$[0,1]$ computed with the Dynamic Programming software differed slightly from one another, while the computed
elastic shape distances between the curves in $3-$dimensional space used as input to the software differed from
one another as well and were not exactly close to zero.
The graph of the optimal diffeomorphism computed the $51^{st}$ time the Dynamic Programming software
was executed is shown in Figure~\ref{F:c312best} together with results of the partial elastic shape registration
of the pair of surfaces.
Once again we note that perhaps as long as the second surfaces
we have chosen for testing the software are perturbed in the same manner in the $r$ direction, it is
likely the graphs of the optimal diffeomorphisms computed with the Dynamic Programming software
will tend to resemble one another regardless of the surfaces involved.
\begin{figure}
\centering
\begin{tabular}{ccc}
\includegraphics[width=0.3\textwidth]{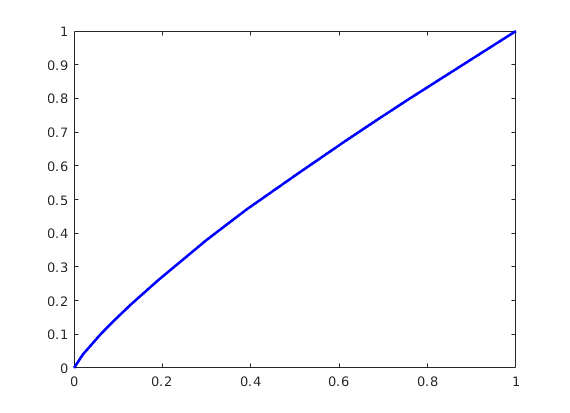}
&
\includegraphics[width=0.3\textwidth]{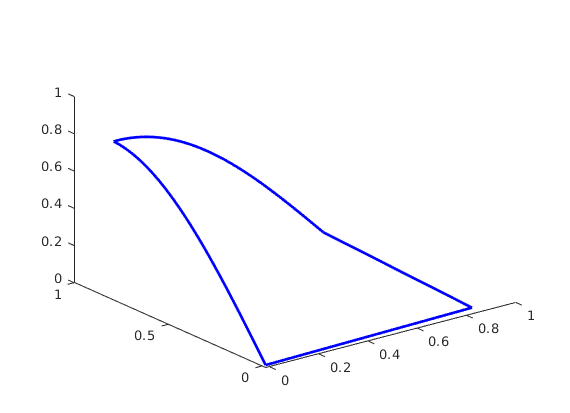}
&
\includegraphics[width=0.3\textwidth]{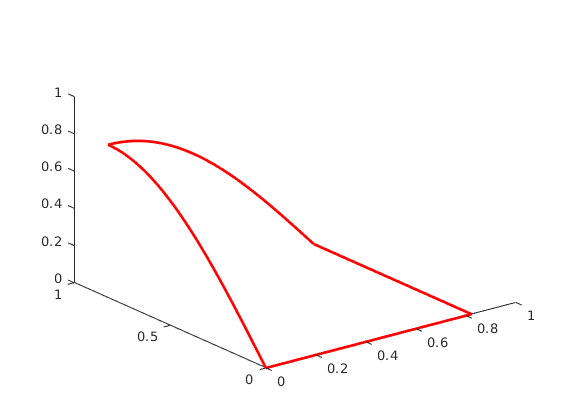}
\end{tabular}
\caption{\label{F:c312best}
For $\gamma(r,t)=(r^{5/4},t^{5/4})$, $(r,t) \in [0,1]\times [0,1]$, views of graph of
optimal diffeomorphism computed the $51^{st}$ time the Dynamic Programming software was executed
(leftmost plot), of boundary of optimally rotated first surface (middle plot), and of optimally
reparametrized second surface (rightmost plot) after computation of partial elastic shape~\mbox{registration}.
}
\end{figure}
\\ \smallskip\\ \noindent
{\bf\large Summary}
\\ \smallskip\\
In this paper we have presented an algorithm for computing, using Dynamic Programming, a partial
elastic shape registration of two simple surfaces in $3-$dimensional space together with the elastic
shape distance between them associated with the partial registration. The algorithm we have presented minimizes
a distance function of the surfaces in terms of rotations of one of the surfaces and a special subset of the set
of reparametrizations of the other surface, the optimization over reparametrizations based on the computation,
using Dynamic Programming, of the elastic shape registration of pairs of simple curves in $3-$dimensional space
contained in the surfaces. This algorithm does not necessarily compute an optimal elastic shape registration
of the surfaces together with the exact elastic shape distance between them, but perhaps a registration and
a distance closer to optimal than those obtained with an algorithm based on a gradient approach over the
entire set of reparametrizations of one of the surfaces. In fact we propose that when computing the elastic
shape registration of two simple surfaces and the elastic shape distance between them with an algorithm based
on a gradient approach for optimizing over the entire set of reparametrizations of one of the surfaces, to use
as the input initial solution the rotation and the reparametrization computed with our proposed algorithm.
Finally, we note, promising results from computations with the implementation of our methods applied on three
simple kinds of $3-$dimensional surfaces, have been presented in this paper. A link to the software package,
etc., has been given as well.

\end{document}